\documentclass[envcountsame]{svmult}
\usepackage{graphicx,amssymb,amsmath}
\smartqed\overfullrule=5pt \long\def\comment#1{}

\let\leq=\leqslant\let\le=\leqslant
\let\geq=\geqslant\let\ge=\geqslant

\begin{document}

\title*{Two recursive decompositions of Brownian bridge
related to the asymptotics of random mappings}

\titlerunning{Two recursive decompositions of Brownian
bridge}

\author{David Aldous and Jim Pitman\thanks{Research
supported in part by N.S.F. Grants DMS-9970901 and
DMS-0071448}}

\institute{Department of Statistics, University of
California\\ 367 Evans Hall \# 3860, Berkeley, CA
94720-3860\\ \texttt{e-mail: aldous@stat.berkeley.edu, pitman@stat.berkeley.edu}}

\maketitle

\begin{abstract} Aldous and Pitman (1994) studied
asymptotic distributions as $n \to \infty$, of various
functionals of a uniform random mapping of the set $\{1,
\ldots, n \}$, by constructing a {\em mapping-walk}\/ and
showing these random walks converge weakly to a reflecting
Brownian bridge.  Two different ways to encode a mapping as
a walk lead to two different decompositions of the Brownian
bridge, each defined by cutting the path of the bridge
at an increasing sequence of recursively defined random
times in the zero set of the bridge.  The random mapping
asymptotics entail some remarkable identities involving
the random occupation measures of the bridge fragments
defined by these decompositions. We derive various
extensions of these identities for Brownian and Bessel
bridges, and characterize the distributions of various
path fragments involved, using the L\'evy--It\^o theory
of Poisson processes of excursions for a self-similar
Markov process whose zero set is the range of a stable
subordinator of index $\alpha \in (0,1)$.

\keywords{Brownian bridge, Brownian excursion, local time,
occupation measure, stable subordinator, self-similar
Markov process, Bessel process, path decomposition,
Poisson--Dirichlet distribution, pseudo-bridge, random
mapping, size-biased sampling, weak convergence,
exchangeable interval partition.} \end{abstract}

\comment{\tableofcontents\newpage}

\section{Introduction} In a previous paper
\cite{aldous:ap92} we showed how features of a uniformly
distributed random mapping $M_n$, from $[n]:= \{1,2,
\ldots, n \}$ to itself, could be encoded as functionals of
a particular non-Markovian random walk on the non-negative
integers.  This {\em mapping-walk}, suitably rescaled,
converges weakly in $C[0,1]$ as $n \to \infty$ to the
distribution of the reflecting Brownian bridge defined by
the absolute value of a {\em standard Brownian bridge}\/
$B^{\rm br}$ with $B^{\rm br}_0 = B^{\rm br}_1 = 0$
obtained by conditioning a standard Brownian motion
$B$ on $B_1 = 0$.  Two important features of a mapping
are the vector of sizes of connected components of its
digraph, and the vector of sizes of cycles in its digraph.
Results of \cite{aldous:ap92} imply that for a uniform
random mapping, as $n \to \infty$, the component sizes
rescaled by $n$, jointly with corresponding cycle sizes
rescaled by $\sqrt{n}$, converge in distribution to
a limiting bivariate sequence of random variables
\comment{\begin{equation}\label{aldous:bivseq}}
$(\lambda_{I_j}, L^0_{I_j})_{j = 1,2, \ldots}$
\comment{\end{equation}} where $(I_j)_{j = 1,2,
\ldots}$ is a random interval partition of $[0,1]$, with
$\lambda_{I_j}$ the length of $I_j$ and $L^0_{I_j}$
the increment of local time of $B^{\rm br}$ at $0$ over
the interval $I_j$.  With the convention for ordering
connected components of the mapping digraph used in
\cite{aldous:ap92}, the limiting interval partition is
$(I_j) = (I^D_j)$, according to the following definition.
Here, and throughout the paper, $U$, $U_1$, $U_2$, $\ldots$
denotes a sequence of independent uniform $(0,1)$
variables, independent of $B^{\rm br}$.

\begin{definition}[the\/ $D$-partition
\cite{aldous:ap92}]\label{aldous:ddef} Let\/ $I^D_j :=
[D_{V_{j-1}}, D_{V_j}  ]$ where\/ $V_0 = D_{V_0}= 0$ and\/
$V_j$ is defined inductively along with the \/ $D_{V_j}$
for\/ $j \ge 1$ as follows: given that\/ $D_{V_i}$ and\/
$V_i$ have been defined for\/ $0 \le i < j$, let
$$
	V_j := D_{V_{j-1}} + U_j(1 - D_{V_{j-1}}),
$$
so\/ $V_j$
is uniform on\/ $[D_{V_{j-1}},1]$ given\/ $B^{\rm br}$
and\/ $(V_i, D_{V_i})$ for\/ $0 \le i < j$, and let
$$
	D_{V_j}:= \inf \bigl\{t \ge V_j : B^{\rm br}_t = 0\bigr\}. 
$$
\end{definition}

On the other hand, a variation of the main result of
\cite{aldous:ap92} shows that with a different ordering
convention, the mapping component sizes rescaled by $n$,
jointly with their cycle sizes rescaled by $\sqrt{n}$,
have a limit distribution specified by the sequence of
lengths and Brownian local times $(\lambda_{I_j}, L^0_
{I_j})_{j = 1,2, \ldots}$ a differently defined limiting
interval partition.  This is the partition $(I_j) =
(I^T_j)$ defined as follows using the local time process
$(L^0_u, 0 \le u \le 1)$ of $B^{\rm br}$ at $0$:

\begin{definition}[the\/ $T$-partition]\label{aldous:tdef}
Let\/ $I^T_j :=  [T_{j-1}, T_j  ]$ where\/ $T_0:= 0$,
$\widehat{V}_0:= 0$, and for\/ $j \ge 1$
\begin{equation}
\label{aldous:tkdef}
	\widehat{V}_j := 1 - \prod_{i =1}^j (1-U_i),
\end{equation}
so\/ $\widehat{V}_j$ is
uniform on\/ $[\widehat{V}_{j-1},1]$ given\/ $B^{\rm br}$
and\/ $(\widehat{V}_i,T_i)$ for\/ $0 \le i < j$, and
$$
	T_j := \inf \bigl\{u : L^0_u /L^0_1 > \widehat{V}_j \bigr\}. 
$$
\end{definition}

For each of these two random interval partitions $(I_j)$
we are interested in the distribution of the bivariate
sequence of lengths and local times $\smash{(\lambda_{I_j},
L^0_{I_j})_{j = 1,2, \ldots}}$ and the distribution
of the associated path fragments $B^{\rm br}[I_j]$ and
standardized fragments $B^{\rm br}_*[I_j]$.  Here for a
process $X:= (X_t, t \in J)$ parameterized by an interval
$J$, and $I = [G_I,D_I]$ a subinterval of $J$ with length
$\lambda_{I} := D_I - G_I > 0$, we denote by $X[I]$ or
$X[G_I,D_I]$ the {\em fragment of $X$  on $I$}, that is
the process
\begin{equation}
\label{aldous:shifty}
	X[I]_u := X_{G_I + u} \qquad (0 \le u \le \lambda_{I}).
\end{equation}
We denote by $X_* [I]$ or $X_* [G_I,D_I]$
the {\em standardized fragment of\/ $X$ on\/ $I$}, defined
by the {\em Brownian scaling operation}
\begin{equation}
\label{aldous:bscale}
	X_* [I]_u := \frac{X [I]_{u\lambda_{I}}}{\sqrt{\lambda_{I}}}
	:= \frac{X_{G_I+u \lambda_{I}}}{\sqrt{\lambda_{I}}}
	\qquad (0 \leq u\leq 1). 
\end{equation}

Figure 1 illustrates these definitions for a typical
path of $X = B^{\rm br}$.  Note that the first
interval $I^D_1$ of the $D$-partition ends at the
time $D_{U_1}$ of the first zero of $B^{\rm br}$ after
a uniform$(0,1)$-distributed time $U_1$, whereas the
first interval $I^T_1$ of the $T$-partition ends at the
time $T_1$ when the local time of $B^{\rm br}$ at $0$
has reached a uniform$(0,1)$-distributed fraction of
its ultimate value.  As illustrated in Figure 1, the
associated fragments of $B^{\rm br}$ are qualitatively
different: $B^{\rm br}_*[I^D_1]$ ends with an excursion
while $B^{\rm br}_*[I^T_1]$ does not.  \comment{xxx Figure
1 needs some polishing; hats on $\hat{V}_i$; mark the
$D_i$; mark unit interval $[0,1]$ somewhere; play with
line thickness. Change notation to $L^0_t$.}

$$
\includegraphics[bb=80 85 575 625,width=0.99\textwidth]{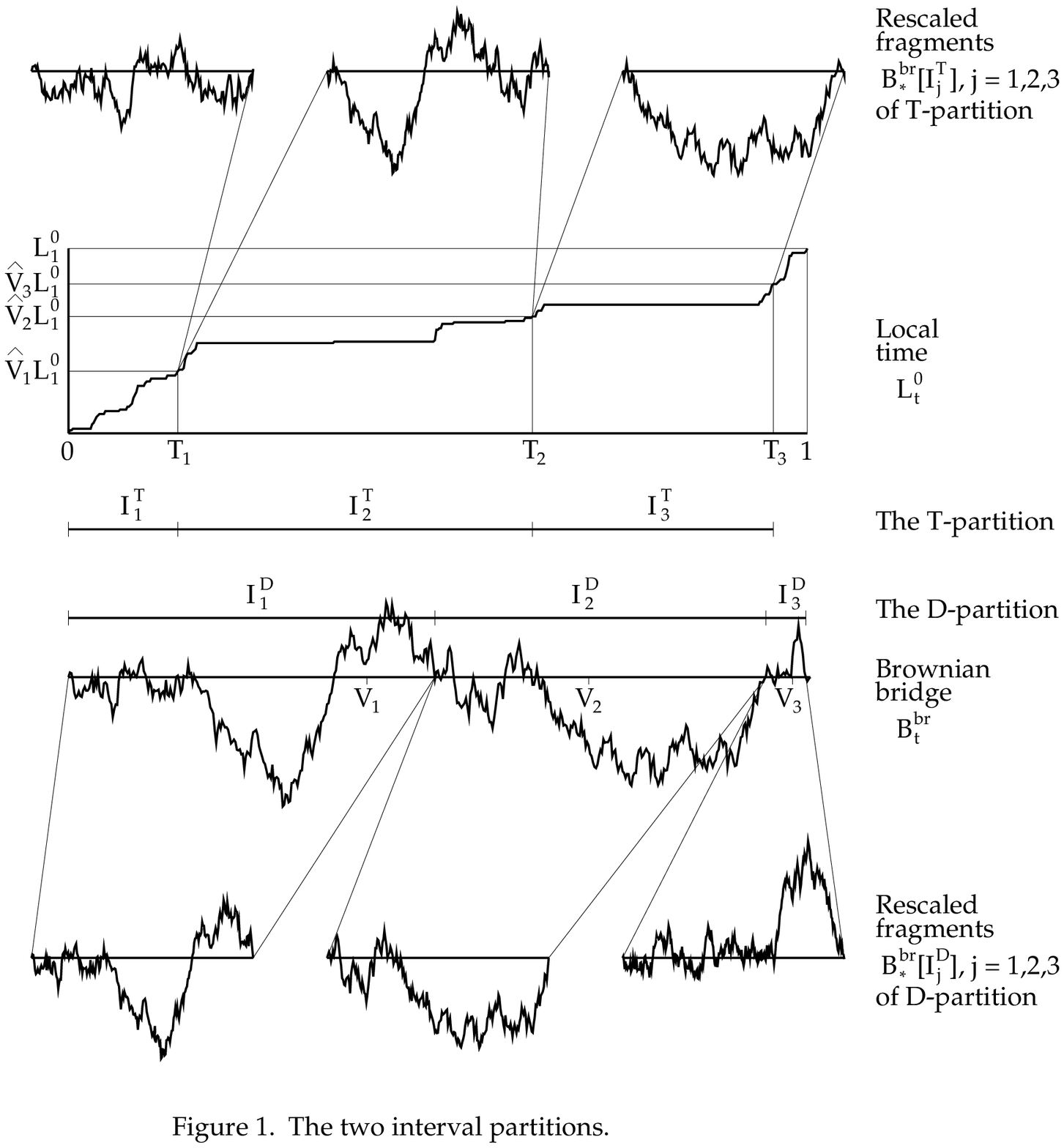}
$$

Despite this difference between the fragments of $B^{\rm
br}$ over the $D$- and $T$-partitions, the random
mapping asymptotics have the following corollary.  Let
$(I^D_{(j)})$ and $(I^T_{(j)})$ denote the length-ranked
$D$-partition and the length-ranked $T$-partition
respectively, meaning $I^D_{(j)}$ is the $j$th longest
interval in the $D$-partition, and $I^T_{(j)}$ is the
$j$th longest interval in the $T$-partition.

\begin{theorem}\label{aldous:crl1} Considering the four
bivariate sequences\/ $(\lambda_{I_j}, L^0_{I_j})_{j =
1,2, \ldots}$ of lengths and bridge local times at\/ $0$,
for\/ $(I_j)$ one of the four random interval partitions
of\/ $[0,1]$ defined by\/ $(I^D_j)$, $(I^D_{(j)})$,
$(I^T_j)$ or\/ $(I^T_{(j)})$,

\noindent {\em (i)}\enspace the bivariate sequence
has the same distribution for\/ $(I^D_{(j)})$ as for\/
$(I^T_{(j)})$;

\noindent {\em (ii)}\enspace the bivariate sequence for\/
$(I^D_{j})$ is the bivariate sequence for\/ $(I^D_{(j)})$
in a length-biased order;

\noindent {\em (iii)}\enspace the bivariate sequence for\/
$(I^T_{j})$ is the bivariate sequence for\/ $(I^T_{(j)})$
in an\/ $L^0$-biased order;

\noindent {\em (iv)}\enspace the sequence of local times\/
$(L^0_{I_j})$ has the same distribution for\/ $(I^D_{j})$
as for\/ $(I^T_{j})$, whereas the sequence of lengths\/
$(\lambda_{I_j})$ does not.  \end{theorem}

See \cite{aldous:ppy92,aldous:jp96bl} for background about
size-biased random orderings.  To illustrate the meaning
of (iii) for instance, for each $k \ge 1$, conditionally
given the entire bivariate sequence $(\lambda_{I^T_{(j)}},
L^0_{I^T_{(j)}})_{j =1,2, \ldots}$, the probability
of the event $(I^T_1 = I^T_{(k)})$ is $L^0_{I^T_{(k)}}
/ L^0_1 $, where \comment{$$
P(I^T_{1} = I^T_{(k)}
\,|\,(\lambda_{I^T_{(j)}}, L^0_{I^T_{(j)}})_{j=1,2,
\ldots}) = \frac{L^0_{I^T_{(k)}}}{L^0_1}
$$} $L^0_1 =
\sum_j L^0_{I^T_{(j)}} $ almost surely. And given
also $(I^T_1 = I^T_{(k)})$, for each $m \ge 1$ with
$m \ne k$ the probability of the event $(I^T_{2}
= I^T_{(m)})$ is \comment{$$
P(I^T_{2} = I^T_{(m)}
\,|\,(\lambda_{I^T_{(j)}}, L^0_{I^T_{(j)}})_{j=1,2,
\ldots}, I^T_{1} = I^T_{(k)}) =
\frac{L^0_{I^T_{(m)}}}{L^0_1 -  L^0_{I^T_{(k)}}}
$$} $L^0_{I^T_{(m)}} /
(L^0_1 -  L^0_{I^T_{1}})$ and so on.  Put another way,
parts (i)-(iii) of the corollary state that the bivariate
sequence $(\lambda_{I_j}, L^0_{I_j})_{j = 1,2, \ldots}$
for $(I_j)= (I^D_j)$ is distributed like a length-biased
rearrangement of the bivariate sequence for $(I_j) =
(I^T_j)$, which is in turn distributed like an $L^0$-biased
rearrangement of the bivariate sequence for $(I_j) =
(I^D_j)$.  Consequently, the distribution of any one of
the four bivariate sequences determines the distribution
of each of the others.

The rest of this paper is organized as follows.
Section~\ref{aldous:sec-RM} explains how we discovered
Theorem~\ref{aldous:crl1} by consideration of random
mapping asymptotics. We recall the theorem from
\cite{aldous:ap92} which describes the asymptotics of
mapping-walks in terms of the fragments of $B^{\rm br}$
defined by the $D$-partition, and present the companion
result, for a different ordering of components, where the
limit involves the fragments of $B^{\rm br}$ defined by
the $T$-partition.  Section~\ref{aldous:sec-state} lays
out our results regarding the decomposition of $B^{\rm
br}$ into path fragments associated with the $D$- and
$T$-partitions, in a way which does not depend on the
random mapping asymptotics.  In particular, we describe
the three different distributions of bivariate sequences
featuring in the three parts of Theorem~\ref{aldous:crl1}.
We formulate and prove these results more generally,
for $B^{\rm br}$ the standardized bridge of a recurrent
self-similar Markov process $B$ whose inverse local
time process at $0$ is a stable subordinator of index
$\alpha$ for some $\alpha \in (0,1)$.  So $\alpha =
1/2$ for $B$ a standard Brownian motion as supposed
in previous paragraphs, and $\alpha = 1 - \delta/2$
for $B$ a Bessel process of dimension $\delta \in (0,2)$.
Some of the results in Section~\ref{aldous:sec-state}, like
Theorem~\ref{aldous:crl1}, can be viewed in the Brownian
case as asymptotic counterparts (under weak convergence of
mapping-walks) of some combinatorial symmetries of random
mappings, discussed in Section~\ref{aldous:comfacts}.
Other results in the Brownian case, especially those
involving the method of Poissonization by random scaling
\cite{aldous:py97max,aldous:py97rh}, are not obvious
from the combinatorial perspective, but provide explicit
limit distributions for functionals of uniform random
mappings.  See also \cite{aldous:ap02d} where we apply
this method to characterize the asymptotic distribution
of the diameter of the digraph of a uniform mapping.
Sections~\ref{aldous:sec-Dpart} and~\ref{aldous:sec-Tpart}
provide some proofs and further details of the
main results in Section~\ref{aldous:sec-state}, while
Section~\ref{aldous:sec-COMP} contains various complements.
In particular, we show in Section~\ref{aldous:sec-IP}
that Theorem~\ref{aldous:crl1} holds even more generally
for interval partitions $(I^D_j)$ and $(I^T_j)$ defined
as before, but with the random zero set of $B^{\rm br}$
replaced by the complement of $\bigcup_j I^{\rm ex}_j$, where
$I^{\rm ex}_j$ is any exchangeable random partition of
$[0,1]$ into an infinite number of intervals, and $(L^0_u,
0 \le u \le 1)$ is the associated local time process,
as defined by Kallenberg \cite{aldous:kal83l}.  This is
the limiting case of a corresponding result for a finite
exchangeable interval partition of $[0,1]$, which we prove
by a combinatorial argument.

In companion papers \cite{aldous:ap01a} and
\cite{aldous:me102} we show that Brownian bridge
asymptotics apply for models of random mappings more
general than the uniform model, in particular for the {\em
p-mapping}\/ model \cite{aldous:op00,aldous:jp01hur}, and
that proofs can be simplified by use of Joyal's bijection
between mappings and trees.  See also \cite{aldous:csp}
for a recent review of the applications of Brownian motion
and Poisson processes to the asymptotics of various kinds
of large combinatorial objects, including partitions,
trees, graphs, permutations, and mappings.

\section{Random Mappings} \label{aldous:sec-RM} In
this section we explain how study of random mappings
led us to consideration of the two interval partitions
of Brownian bridge, and show how the distributions of
path fragments of the bridge defined by these partitions
encode various asymptotic distributions for mappings.
\comment{As mentioned earlier, our focus in this paper
is to prove all those results by process methods, so this
section serves largely as motivation.}

\subsection{Mapping-walks and the two orderings}
\label{aldous:sec-MW} A mapping $M_n :[n] \to [n]$ can
be identified with its digraph of edges $\{(i,M_n (i)),
\, i \in [n]\}$.  The connection between random mappings
and Brownian bridge developed in \cite{aldous:ap92}
can be summarized as follows.  \begin{description}[$\bullet$]

\item[$\bullet$]
A mapping digraph can be decomposed as a collection
of {\em rooted trees}\/ together with extra structure
({\em cycles, basins of attraction}). 

\item[$\bullet$] A rooted
tree can be coded as a discrete {\em tree-walk}, a
walk excursion starting and ending at $0$. 

\item[$\bullet$] Given
some ordering of tree-components, one can concatenate
walk-excur\-sions to define a discrete {\em mapping-walk}\/
which codes $M_n$. 

\item[$\bullet$] For a uniform random mapping,
the induced distribution on tree-compo\-nents is such that
the tree-walks, suitably normalized, converge to Brownian
excursion as the tree size increases to infinity. 

\item[$\bullet$] So
for a uniform random mapping, we expect the mapping-walks,
suitably normalized, to converge to a limit process defined
by some concatenation of Brownian excursions.

\item[$\bullet$] With appropriate choice of ordering, the
limit process is in fact reflecting Brownian bridge.
\end{description}

We now amplify this summary, emphasizing the only subtle
issue -- the choice of ordering.  Fix a mapping $M_n$.
It has a set of {\em cyclic points}
$$
	\mathcal{C}_n:= \bigl\{
	i\in [n]: M_n^k(i) = i \ \mbox{for some}\ k \ge 1\bigr\},
$$
where
$M_n^k$ is the $k$th iterate of $M_n$. \comment{$(M_n(i), i
\in [n])$,} Let $\mathcal{T}_{n,c}$ be the set of vertices
of the (perhaps trivial) tree component of the digraph
with root $c \in \mathcal{C}_n$.  The tree components are
bundled by the disjoint cycles $\mathcal{C}_{n,j} \subseteq
\mathcal{C}_n$ to form the {\em basins of attraction}\/
(connected components) of the mapping digraph, say
\begin{equation}
\label{aldous:basins}
	\mathcal{B}_{n,j}:=
	\bigcup_{c \in \mathcal{C}_{n,j}} \mathcal{T}_{n,c}
	\supseteq \mathcal{C}_{n,j} \quad\mbox{with}\quad
	\bigcup_j \mathcal{B}_{n,j} = [n] \quad\mbox{and}\quad
	\bigcup_j \mathcal{C}_{n,j} = \mathcal{C}_n
\end{equation}
where all three unions are disjoint unions, and the
$\mathcal{B}_{n,j}$ and $\mathcal{C}_{n,j}$ are indexed in
some way by $j = 1, \ldots, K_n$ say.  The construction in
\cite{aldous:ap92} encodes the restriction of the digraph
of $M_n$ to each tree component $\mathcal{T}_{n,c}$ of size
$k$ (that is, with $k$ vertices) by $2k$ steps of a {\em
tree-walk}\/ with increments $\pm1$ on the non-negative
integers. The tree-walk proceeds by a suitable search of
the set $\mathcal{T}_{n,c}$, making an excursion which
starts at $0$ and returns to $0$ for the first time after
$2k$ steps, after reaching a maximum level $1 + h_n(c)$,
where $h_n(c)$ is the maximal height above $c$ of all
vertices of the tree $\mathcal{T}_{n,c}$ with root $c$,
that is
\begin{equation}
\label{aldous:heights}
	h_n(c)
	= \max \bigl\{h : \exists i \!\in\! [n]
	\ \mbox{with}\ M_n^h(i)
	= c \ \mbox{and}\ M_n^{j}(i) \notin \mathcal{C}_n
	\ \mbox{for}\  0 \leq j < h\bigr\}. 
\end{equation}
It was
shown in \cite{aldous:me56} that as $k \to \infty$, the
distribution of the tree-walk for a $k$-vertex random
tree, of the kind contained in the digraph of the uniform
random mapping $M_n$ for $k \le n$, when scaled to have
$2k$ steps of $\pm 1/\sqrt{k}$ per unit time, converges
to the distribution $2 B^{\rm ex}$ for $B^{\rm ex}$ a
standard Brownian excursion.  \comment{whose definition is
recalled in Section~\ref{aldous:sec-BS}.} Subsequent work
\cite{aldous:mm01z} shows that the same result holds for
a variety of codings of trees as walks. Consequently, any
of these codings would serve our purpose in the following
definitions.

We now define a mapping-walk (to code $M_n$) as a
concatenation of its tree-walks, to make a walk of
$2n$ steps starting and ending at $0$ with exactly
$|\mathcal{C}_n|$ returns to $0$, one for each tree
component of the mapping digraph. \comment{then the
distribution of the scaled mapping-walk} To retain useful
information about $M_n$ in the mapping-walk, we want the
definition of the walk to respect the cycle and basin
structure of the mapping.  Here are two orderings that
do so.

\begin{definition}[cycles-first ordering]\label{aldous:cfo}
Fix a mapping\/ $M_n$ from\/ $[n]$ to\/ $[n]$.  If\/ $M_n$
has\/ $K_n$ cycles, first put the cycles in increasing
order of their least elements, say\/ $c_{n,1} < c_{n,2}
< \ldots < c_{n,K_n}$.  Let\/ $\mathcal{C}_{n,j}$ be the
cycle containing\/ $c_{n,j}$, and let\/ $\mathcal{B}_{n,j}$
be the basin containing\/ $\mathcal{C}_{n,j}$.
Within cycles, list the trees around the cycles,
as follows.  If the action of\/ $M_n$ takes\/ $c_{n,j}
\to c_{n,j,1} \to \cdots \to c_{n,j}$ for each\/ $1 \le
j \le K_n$, the tree components\/ $\mathcal{T}_{n,c}$
are listed with\/ $c$ in the order
\begin{equation}
\label{aldous:cforder}
	(\overbrace{c_{n,1,1},\ldots,c_{n,1}}
	^{\mbox{$\mathcal{C}_{n,1}$}},
	\overbrace{c_{n,2,1},\ldots,c_{n,2}}
	^{\mbox{$\mathcal{C}_{n,2}$}},
	\ldots ,
	\overbrace{c_{n,K_n,1},\ldots,c_{n,K_n}}
	^{\mbox{$\mathcal{C}_{n,K_n}$}}). 
\end{equation}
The {\em cycles-first mapping-walk}\/ is obtained by
concatenating the tree walks derived from\/ $M_n$ in this
order.  The {\em cycles-first search of\/ $[n]$} is the
permutation\/ $\sigma: [n] \to [n]$ where\/ $\sigma_j$
is the\/ $j$th vertex of the digraph of\/ $M_n$ which
is visited in the corresponding concatenation of tree
searches.  \end{definition}

\begin{definition}[basins-first ordering
\cite{aldous:ap92}]\label{aldous:bfo} If\/ $M_n$ has\/
$K_n$ cycles, first put the basins\/ $\mathcal{B}_{n,j}$
in increasing order of their least elements, say\/ $1=
b_{n,1} < b_{n,2} < \ldots <b_{n,K_n}$; let\/ $c_{n,j} \in
\mathcal{C}_{n,j}$ be the cyclic point at the root of the
tree component containing\/ $b_{n,j}$. Now list the trees
around the cycles, just as in~\eqref{aldous:cforder},
but for the newly defined\/ $c_{n,j}$ and\/ $c_{n,j,i}$.
Call the corresponding mapping-walk and search of\/
$[n]$ the\/ {\em basins-first mapping-walk} and\/ {\em
basins-first search}.  \end{definition}

Be aware that the meaning of $\mathcal{B}_{n,j}$ and
$\mathcal{C}_{n,j}$ now depends on the ordering convention.
Rather than introduce two separate notations for the two
orderings, we use the same notation for both, and indicate
nearby which ordering is meant.  Whichever ordering, the
definitions of $\mathcal{B}_{n,j}$ and $\mathcal{C}_{n,j}$
are always linked by $\mathcal{B}_{n,j} \supseteq
\mathcal{C}_{n,j}$, and~\eqref{aldous:basins} holds.

Let us briefly observe some similarities between the
two mapping-walks.  For each given basin $B$ of $M_n$
with say $b$ elements, the restriction of $M_n$ to $B$
is encoded in a segment of each walk which equals at $0$
at some time, and returns again to $0$ after $2b$ more
steps.  If the basin contains exactly $c$ cyclic points,
this walk segment of $2b$ steps will be a concatenation of
$c$ excursions away from $0$.  Exactly where this segment
of $2b$ steps appears in the mapping-walk depends on the
ordering convention, as does the ordering of excursions
away from $0$ within the segment of $2b$ steps.  However,
many features of the action of $M_n$ on the basin $B$
are encoded in the same way in the two different stretches
of length $2b$ in the two walks, despite the permutation
of excursions.  One example is the number of elements in
the basin whose height above the cycles is $h$, which
is encoded in either walk as the number of upcrossings
from $h$ to $h+1$ in the stretch of walk of length $2b$
corresponding to that basin.

\subsection{Symmetry properties of random mappings}
\label{aldous:comfacts} We now apply the definitions above
to a uniform random mapping $M_n$.  Of course, the random
partition $\{\mathcal{B}_{n,j} \}_{j = 1, \ldots, K_n}$
of $[n]$, and the random partition $\{\mathcal{C}_{n,j}
\}_{j = 1, \ldots, K_n}$ of $\mathcal{C}_{n}$, are the same
no matter which ordering convention is used.  \comment{One
final trap for the unwary.} Each random partition is {\em
exchangeable}, meaning its distribution is invariant under
the action of a permutation of $[n]$.  Let us spell out
some further symmetry properties, each of which turns out
to have some analog in the limiting Brownian scheme.

\noindent {\bf (a)} The cycles-first ordering\ has the
following very strong symmetry property: conditionally
given $|\mathcal{C}_n| = m$ the tree components in
cycles-first ordering\ form an exchangeable sequence of
$m$ random subsets of $[n]$; moreover this exchangeable
sequence is independent of the sequence of cycle sizes
$|\mathcal{C}_{n,j}|$ with $\sum_{j} |\mathcal{C}_{n,j}|
= m$.  Consequently, given $|\mathcal{C}_n| = m$,
the cycles-first mapping-walk is a concatenation of
$m$ exchangeable excursions away from $0$, and this
mapping-walk is independent of $|\mathcal{C}_{n,j}|, j =
1,2, \ldots, K_n$.

\noindent {\bf (b)} The basins-first ordering\ does not
share the symmetry property above. But it has a different
one: given that the basin $\mathcal{B}_{n,1}$ containing
$1$ has size $|\mathcal{B}_{n,1}| = b$, the action of
$M_n$ on $[n] - \mathcal{B}_{n,1}$ is that of a uniform
random mapping of a set of $n-b$ elements.  So given
$|\mathcal{B}_{n,1}| = b$, the basins-first mapping-walk
decomposes after $2b$ steps into two independent segments:
the first $2b$ steps are distributed like the basins-first
walk for a uniform mapping of $[b]$ conditioned to have a
single basin, and the remaining $2(n-b)$ steps distributed
like the basins-first walk associated with a uniform
mapping of $[n-b]$.

\noindent {\bf (c)} The sequence of basin sizes
$(|\mathcal{B}_{n,j}|, 1 \le j \le K_n)$ does not have the
same distribution for both orderings.  For instance, if
$|\mathcal{B}_{n,1}| = 1$ in the basins-first ordering\
then $|\mathcal{B}_{n,1}| = 1$ in the cycles-first
ordering, but (for $n \geq 3$) not conversely.  So the
distribution of $|\mathcal{B}_{n,1}|$ must be different
in the two orderings.

\noindent {\bf (d)} For a given mapping $M_n$, the sequence
of cycle sizes $(|\mathcal{C}_{n,j}|, 1 \le j \le K_n)$
may be different for the two different orderings. But for
$M_n$ with uniform distribution on $[n]^{[n]}$, the two
sequences of cycle sizes have the same distribution: given
$| \mathcal{C}_n | = m$, either sequence is distributed
like the sizes of cycles of a uniform random permutation of
$[m]$ in the (size-biased) order of least elements of the
cycles.  That is to say, given $| \mathcal{C}_n | = m$, the
distribution of $|\mathcal{C}_{n,1}|$ is uniform on $[m]$;
given $| \mathcal{C}_n | = m$ and $|\mathcal{C}_{n,1}|$
with $| \mathcal{C}_n | - |\mathcal{C}_{n,1}| = m_1$,
the distribution of $|\mathcal{C}_{n,2}|$ is uniform on
$[m_1]$, and so on.  This is a well known property of
uniform random permutations for the cycles-first ordering,
and was shown for the basins-first ordering in \cite[Lemma
22]{aldous:ap92}.

\subsection{Brownian asymptotics for the mapping-walks}
\label{aldous:asym} We now come to the main point of
Section~\ref{aldous:sec-RM}: the definitions of the
interval partitions of Brownian bridge are motivated by
the following theorem.

\begin{theorem}\label{aldous:thconv} The scaled
mapping-walk\/ $({M_{u}^{[n]}}, 0 \le u \le 1)$, with\/
$2n$ steps of\/ $\pm 1/\sqrt{n}$ per unit time, for
either the cycles-first or the basins-first ordering\ of
excursions corresponding to tree components, converges
in distribution to\/ $2 |B^{\rm br}|$ jointly with
\begin{equation}
\label{aldous:cyconv}
	\frac{|\mathcal{C}_{n}|}{\sqrt{n}}
	\stackrel{\!d}{\longrightarrow}  L^0_1
\end{equation}
where\/ $(L^0_u, 0 \le u \le 1)$ is
the process of local time at\/ $0$ of\/ $B^{\rm br}$,
normalized so that\/ $P(L^0_1 > \ell) = \E^{-\ell^2/2}$.
Moreover,

\noindent {\em (i)}\enspace for
the cycles-first ordering, with the cycles\/
$\mathcal{B}_{n,j}$ in order of their least elements,
these two limits in distribution hold jointly with
\begin{equation}
\label{aldous:wconv}
	\biggl(\frac{|\mathcal{B}_{n,j}|}{n}\,,\,
	\frac{|\mathcal{C}_{n,j}|}{\sqrt{n}} \biggr)
	\stackrel{\!d}{\longrightarrow}
	\bigl(\lambda_{I_j}, L^0_{I_j}\bigr)
\end{equation}
as\/ $j$ varies, where the limits are the
lengths and increments of local time of\/ $B^{\rm br}$
at\/ $0$ associated with the interval partition\/ $(I_j) :=
(I^T_j)$; whereas

\noindent {\em (ii) \cite{aldous:ap92}}\enspace
for the basins-first ordering, with the basins\/
$\mathcal{B}_{n,j}$ listed in order of their least
elements, the same is true, provided the limiting interval
partition is defined instead by\/ $(I_j) := (I^D_j)$.
\end{theorem}

The result for basins-first ordering\  is part
of \cite[Theorem 8]{aldous:ap92}.  The variant for
cycles-first ordering\ can be established by a variation
of the argument in \cite{aldous:ap92}, exploiting
the exchangeability property of the cycles-first
ordering (Section~\ref{aldous:comfacts}(a) instead
of Section~\ref{aldous:comfacts} (b)).  See also
\cite{aldous:biane94bb} and \cite{aldous:me102}
for alternate approaches to the basic result of
\cite{aldous:ap92}.

\comment{xxx rest of this section is quite chaotic. How to
improve? Remark (e) seems most important and should come
earlier. I would like to see a genuine sketch of proof
for the CFO, plus an explanation of how the result for
BFO could be deduced from the CFO result and our Brownian
results, without further combinatorial considerations.}

We now explain how we first discovered some
of the facts about Brownian bridge presented
in Theorem~\ref{aldous:crl1} by consideration of
Theorem~\ref{aldous:thconv} and the symmetry properties
of Section~\ref{aldous:comfacts}.  The arguments below
are not part of the formal development in this paper.
Indeed we show in Section~\ref{aldous:sec-COMP} that
the results of Theorem~\ref{aldous:crl1} hold much more
generally, so these results do not really involve much
of the rich combinatorial structure of mapping digraphs
involved in Theorem~\ref{aldous:thconv}.

\noindent {\bf (a)} In the basins-first ordering, the first
basin is by definition the basin containing element $1$,
and its walk-segment ends at the first time that the
walk returns to $0$ after the basins-first search has
reached element $1$.  Suppose we could replace element
$1$ by a uniform random element, so the walk-segment
corresponds asymptotically to the walk-segment ending
at the first time of reaching $0$ after a uniform random
time on $[0,2n]$.  Rescaling, this corresponds to the time
interval $[0,D_{V_1}]$ in Definition~\ref{aldous:ddef}.
Of course it is not obvious, and indeed is somewhat
counter-intuitive, that replacing element $1$ by a uniform
random element will preserve length of walk-segment.  But
this is what eventually will emerge from our calculations.

Now consider the cycles-first ordering.  The first basin is
by definition the basin containing the smallest-numbered
cyclic element $c_{n,1}$, and its walk-segment ends
at the first time after reaching element $c_{n,1}$
that the walk returns to $0$.  Suppose as before (and
again this is not obvious) one can replace element
$c_{n,1}$ by a uniform random {\em cyclic}\/ element,
so the walk-segment corresponds asymptotically to the
walk-segment ending at the first time of reaching $0$
after \comment{a uniform random time on $[0,2n]$.}
visiting $U^* |\mathcal{C}_n|$ cyclic vertices, where
$U^*$ has uniform$[0,1]$ distribution.  \comment{\tt xxx
I dont understand this argument either. I would prefer to
refer explicitly to Section~\ref{aldous:comfacts}(d), plus
discussion of convergence of the process counting mapping
walk visits to $0$ to the bridge local time process, which
should perhaps be added to Theorem~\ref{aldous:thconv}. }
Rescaling, this corresponds to the time interval $[0,T_1]$
in Definition~\ref{aldous:tdef}.

\noindent {\bf (b)} The recursive property of the
basins-first ordering\ in Section~\ref{aldous:comfacts}(b)
plainly corresponds, under the asymptotics of
Theorem~\ref{aldous:thconv}, to the recursive decomposition
of Brownian bridge at time $D_{V_1}$ described later in
Lemma~\ref{aldous:lmm-Drec}.

\noindent {\bf (c)} In Section~\ref{aldous:comfacts}(c)
we observed that the distribution of $\mathcal{B}_{n,1}$
was different in the two orderings.  This difference
persists in the limit: Theorem~\ref{aldous:thconv} and
the calculation below~\eqref{aldous:neq} imply
$$
\lim_n
n^{-1} E | \mathcal{B}_{n,1}| = \left\{\begin{array}{lll}
E(D_{V_1}) = 2/3 && (\mbox{for the basins-first ordering})
\\ E(T_1) = 1/2 && (\mbox{for the cycles-first ordering}).
\end{array} \right. 
$$
\noindent {\bf (d)} It is well
known \cite{aldous:ve77} that the asymptotic distribution
as $n \to \infty$ of the fractions of elements in cycles
of a random permutation of $[n]$, with the cycles in order
of their least elements, (which amounts to a size-biased
random order by exchangeability), is the {\em uniform
stick-breaking sequence}\/ $U_j \prod_{i=1}^{j-1} (1 -
U_i)$. \comment{appearing here is well known to} So the
convergence in distribution \eqref{aldous:cyconv} of
$|\mathcal{C}_n|/\sqrt{n}$ to $L^0_1$, and the ``uniform
random permutation" feature of the cyclic decomposition
(Section~\ref{aldous:comfacts}(d)), combine to show that
with {\em either}\/ ordering $\smash{|\mathcal{C}_{n,j}|/\sqrt{n}
\stackrel{d}{\rightarrow} L^0_{I_j}}$ with the same joint
distribution:
\begin{equation}
\label{aldous:locstick}
	\bigl(L^0_{I_j}, j \geq 1\bigr)
	\stackrel{d}{=}
	\Biggl(L^0_1
	U_j \prod_{i=1}^{j-1} (1-U_i), \ j \geq 1\Biggr)
\end{equation}
for both $I_j = I^D_j$ and $I_j = I^T_j$.
This is part (iv) of Theorem~\ref{aldous:crl1}, which
is generalized later by~\eqref{aldous:Lsame} and
Theorem~\ref{aldous:thmex}.

\noindent {\bf (e)} Let $\mathcal{B}_{n,(j)}$ be the $j$th
largest basin of $M_n$, with some arbitrary convention for
breaking ties, and let $\mathcal{C}_{n,(j)}$ be the cycle
contained in $\mathcal{B}_{n,(j)}$.  It follows immediately
from the convergence in distribution~\eqref{aldous:wconv}
that
\begin{equation}
\label{aldous:wconv1}
	\biggl(\frac{|\mathcal{B}_{n,(j)}|}{n}\,,\,
	\frac{|\mathcal{C}_{n,(j)}|}{\sqrt{n}} \biggr)
	\stackrel{\!d}{\longrightarrow}
	\bigl(\lambda_{I_{(j)}},
	L^0_{I_{(j)}}\bigr)
\end{equation}
jointly as $j$ varies,
where $I_{(j)}$ is the length-ranked interval partition
derived from either $(I^D_j)$ or $(I^T_j)$.  This is part
(i) of Theorem~\ref{aldous:crl1}.  By exchangeability
considerations, before passage to the limit the bivariate
sequence in~\eqref{aldous:wconv} as $j$ varies is that
in~\eqref{aldous:wconv1} biased by cycle-size in the
cycles-first order and biased by basin-size in the
basins-first ordering.  Hence the conclusions of parts
(ii) and (iii) of Theorem \ref{aldous:crl1}, by a
straightforward passage to the limit.

\noindent {\bf (f)} Due to
Section~\ref{aldous:comfacts}(a), it makes no
difference to anything if in the cycles-first
ordering we replace the ordering within the $j$th
cycle $c_{n,j,1}, c_{n,j,2}, \ldots,   c_{n,j}$
by the possibly more natural $c_{n,j}, c_{n,j,1},
c_{n,j,2}, \ldots,   c_{n,j,|\mathcal{C}_{n,j}|-1}$.
But in the basins-first ordering, this innocent looking
change would spoil convergence to $2 |B^{\rm br}|$.
This is because in the basins-first ordering\ the tree
with root $c_{n,1}$ is the tree containing $1$, which
is a size-biased choice from the exchangeable random
partition of $[n]$ into tree components.  As such, it
tends to be a big tree.  \comment{which we know in the
limit consumes on average of about $n/3$ vertices.} In
fact, results from \cite{aldous:ap92} imply that, if the
mapping-walk is started by the excursion coding the tree
rooted at $c_{n,1}$, the limit process will start with a
zero free interval whose length is distributed as $D_U
- G_U$ in Lemma~\ref{aldous:lmmP2} below for $\alpha =
1/2$. Such a process is obviously not $2 |B^{\rm br}|$
or any other familiar Brownian process.

\noindent {\bf (g)} The proof of
Theorem~\ref{aldous:thconv} yields more information about
the asymptotic sizes of tree components than can be deduced
from the statement of that theorem.  For instance, if $|
\mathcal{T}_{n,(i)}|$ are the ranked sizes of the tree
components of $M_n$, and $H_{n,i}$ are the corresponding
maximal tree heights, as in~\eqref{aldous:heights},
then $ (|\mathcal{T}_{n,(i)}|/n, H_{n,i}/\sqrt{n})_{i =
1,2, \cdots}$ converges in distribution to the sequence
of ranked lengths and corresponding maximal heights of
excursions of $2|B^{\rm br}|$, whose distribution was
described in \cite[Theorem 1 and Example 8]{aldous:py97rh}.
If only the tree components of $\mathcal{B}_{n,j}$ were
considered, the limit would be derived from excursions of
$B^{\rm br}$ over the appropriate random interval $I_j$
as in Theorem~\ref{aldous:thconv}, with joint convergence
as $j$ varies.

\section{The bridge decompositions}
\label{aldous:sec-state} This section presents our main
results for the $D$- and $T$-partitions.  For ease of
comparison, the results are presented together here, with
outlines of the proofs.  Some proofs and further details
are deferred to Section~\ref{aldous:sec-Dpart} for the
$D$-partition, and to Section~\ref{aldous:sec-Tpart}
for the $T$-partition.  Our primary interest is the
analysis of the $D$- and $T$-partitions derived from a
standard Brownian bridge, and the connections between
these random partitions and the asymptotics of random
mappings discussed in Section~\ref{aldous:sec-RM}.
But we find that our analysis applies just as well to
the $D$- and $T$-partitions for a standardized bridge
$B^{\rm br}$ derived from $B$ a recurrent self-similar
Markov process whose inverse local time process at $0$
is a stable subordinator of index $\alpha$ for some
$\alpha \in (0,1)$.  Readers who don't care about this
generalization can assume throughout this section that $B$
is standard one-dimensional Brownian motion, and \hbox{$\alpha =
\beta = 1/2$.}

\subsection{General framework} Following Pitman--Yor
\cite[\S 2]{aldous:py97rh}, we make the following basic
assumptions:

\begin{description}[$\bullet$]

\item[$\bullet$] $B:= (B_t, t \ge 0)$ is a real or
vector-valued strong Markov process, started at $B_0 =
0$, with state space a cone contained in $\mathbb{R}^d$
for some $d = 1,2, \ldots$, and c\`adl\`ag paths.

\item[$\bullet$] $B$ is $\beta$-self-similar for some real
$\beta$. That is to say, if $B_* [0,t]$ now denotes
the standardized process derived from $B$ on $[0,t]$
as in~\eqref{aldous:bscale}, using $\lambda_{I}^\beta$
instead of $\sqrt{\lambda_{I}}$ in the denominator, then
$\smash{B_* [0,t] \stackrel{d}{=} B[0,1]}$ for all $t >0$.

\item[$\bullet$] The point $0$ is a regular recurrent point for
$B$, meaning that almost surely both $0$ and $\infty$
are points of accumulation of the zero set of $B$.
\end{description} As a well known consequence of these
assumptions \cite{aldous:gp80,aldous:py97rh}, there exists
a continuous local time process for $B$ at $0$, say $(L_t^0
(B), t \ge 0)$, whose inverse process
$$
	\tau_\ell :=
	\inf \{t : L_{t}^0 (B) > \ell \} \qquad (\ell \ge 0)
$$
is a stable subordinator of index $\alpha$ for some $\alpha
\in (0,1)$.  That is
\begin{equation}
\label{aldous:taul}
	E \exp (- \xi \tau_\ell) = \exp (- \ell c \xi^\alpha)
	\qquad (\xi \ge 0)
\end{equation}
for some $c >0$, in which
case
\begin{equation}
\label{aldous:limitl}
	L_{t}^0 (B)
	= \frac{\Gamma(1-\alpha)}{c} \,\lim_{\varepsilon \to 0}
	\varepsilon^{\alpha} N_{t, \varepsilon} (B)
\end{equation}
uniformly for bounded $t$ almost surely, where $N_{t,
\varepsilon} (B)$ is the number of excursion intervals of
$B$ in $[0,t]$ whose length is greater than $\varepsilon$.
Formula~\eqref{aldous:limitl} can be then used with
$X$ instead of $B$ to define $L_t^0(X)$ for various
other processes $X$ derived from $B$ by conditioning or
scaling, such as the standardized bridge $B^{\rm br}$
introduced in the next paragraph. As a consequence
of~\eqref{aldous:limitl} with $X$ instead of $B$,
there is following basic {\em $\alpha$-scaling rule}\/
for such local time processes: for $I = [G_I,D_I]$
a random subinterval of length $\lambda_{I}:= D_I -
G_I$ contained in the time domain of $X$, and $L_I^{0}
(X):= L_{D_I}^0 (X) - L_{G_I}^0(X)$,
\begin{equation}
\label{aldous:locscale}
	L_I^{0} (X) = \lambda_{I}^\alpha\,L_{1}^{0}  (X_* [I]). 
\end{equation}
Associated
with the self-similar Markov process $B$ are corresponding
distributions of a {\em standard\/ $B$-bridge}\/ $B^{\rm
br}$, a {\em standard\/ $B$-excursion}\/ $B^{\rm ex}$,
and a {\em standard\/ $B$-meander}\/ $B^{\rm me}$, defined
by the following identities in distribution, valid for
all $t >0$:
\begin{equation}
\label{aldous:scale}
	B_*[0,G_t] \stackrel{d}{=} B^{\rm br}; \qquad B_* [G_t,D_t] 
	\stackrel{d}{=}  B^{\rm ex} ; \qquad B_* [G_t, t]
	\stackrel{d}{=} B^{\rm me}
\end{equation}
where $G_t:=
G_t(B)$, $D_t:= D_t(B)$, and for any process $X$ we use
the notation
\begin{eqnarray} \label{aldous:gtdef}
	G_t(X)
	:= & \sup \{u < t : X_u = 0 \} \\ \label{aldous:dtdef}
	D_t (X) := &  \inf \{u \ge t : X_u = 0 \}.
\end{eqnarray}

See \cite{aldous:bp92} for a review of properties of
$B^{\rm br}$, $B^{\rm ex}$ and $B^{\rm me}$ in the {\em
Brownian case}\/ when $B$ is Brownian motion with state
space $\mathbb{R}$, and $\beta  = \alpha = 1/2$.
See \cite[\S 3]{aldous:ppy92} and \cite{aldous:py97max}
for some treatment of $B^{\rm br}$ and $B^{\rm ex}$
in the {\em Bessel case}\/ when $B$ with state space
$\mathbb{R}_{\ge 0}$ is a recurrent Bessel process of
dimension $\delta = 2 - 2 \alpha \in (0,2)$, and $\beta =
1/2$.  Other examples are provided by recurrent stable
L\'evy processes \cite{aldous:bertoin-levy}, symmetrized
or skew Bessel processes \cite{aldous:wat93}, and Walsh
processes \cite{aldous:bpy89a,aldous:bpy89b}.

According to the L\'evy--It\^o theory of excursions
of $B$, applied to the standard $B$-bridge as in
\cite{aldous:py92,aldous:py97rh}, if $(I^{\rm ex}_j)$
is the interval partition of $[0,1]$ defined by the
excursion intervals of $B^{\rm br}$ in length-ranked
order, then the processes $B^{\rm br}_*[I^{\rm ex}_j]$
are i.i.d.\ copies of $B^{\rm ex}$, independent of
$(I^{\rm ex}_j)$, which is an {\em exchangeable interval
partition}\/ in the sense of \cite{aldous:kal83l} recalled
in Section~\ref{aldous:sec-IP}.  Moreover, the distribution
of ranked lengths $(\lambda_{I^{\rm ex}_j})$ depends only
on $\alpha$, as described in \cite[(16)]{aldous:py95pd2}
and \cite[Example 8]{aldous:py97rh}.  This general
excursion decomposition of $B^{\rm br}$ implies that
various results known for Bessel bridges hold also in
the present general setting, and we take this for granted
without further comment.

\subsection{Main Results} All results of this section are
presented with the notation and general framework of the
previous section:  $B^{\rm br}$ is the standard $B$-bridge
derived from a self-similar recurrent Markov process $B$
whose continuous local time process $(L_t^0(B), t \ge 0)$
is the inverse of a stable subordinator $(\tau_\ell,
\ell \ge 0)$ of index $\alpha \in (0,1)$.  The $D$-
and $T$-partitions are defined in terms of $B^{\rm
br}$ and its local time process at $0$, according to
Definitions~\ref{aldous:ddef} and~\ref{aldous:tdef}.

Theorem~\ref{aldous:crl1}, presented in the introduction
in the Brownian case, is true in the more general framework
of this section, as a consequence of the following theorem:

\begin{theorem}\label{aldous:thmDT} For a random interval\/
$I \subseteq [0,1]$, let\/ $\mu_I$ denote the random
occupation measure induced by the path of\/ $B^{\rm
br}[I]$, so for each Borel subset\/ $A$ of the state space
of\/ $B^{\rm br}$
$$
	\mu_I(A) := \int_I 1(B^{\rm br}_t\in A)\,\D t. 
$$

\noindent {\em (i)}\enspace The sequence of occupation
measures\/ $(\mu_{I_j})$ has the same distribution for
each of the two length-ranked partitions\/ $(I_j) =
(I^D_{(j)})$ and\/ $(I_j) = (I^T_{(j)})$.

\noindent {\em (ii)}\enspace For\/ $(I_j) = (I^D_j)$ the
sequence of occupation measures\/ $(\mu_{I_j})$ is in\/
$\lambda_{}$-biased order, where\/ $\lambda_{I_j}$ is the
total mass of the random measure\/ $\mu_{I_j}$.

\noindent {\em (iii)}\enspace For\/ $(I_j) = (I^T_j)$
the sequence of occupation measures\/ $(\mu_{I_j})$ is
in\/ $L^0$-biased order, where\/ $L^0_{I_j} := L^0_{I_j}
(B^{\rm br})$.

\noindent {\em (iv)} For\/ $(I_j)$ equal to any one of the
four interval partitions\/ $(I^D_j)$, $(I^D_{(j)})$,
$(I^T_j)$ or\/ $(I^T_{(j)})$, conditionally given\/
$\lambda_{I_j} = \lambda_j$ and $L^0_{I_j} = \ell_j$ for
all\/ $j = 1,2, \ldots$, the random occupation measures\/
$\mu_{I_j}$, $j = 1,2, \ldots$ are independent, with\/
$\mu_{I_j}$ distributed like the random occupation measure
of a process with the common conditional distribution of
\begin{equation}
\label{aldous:switch}
	(B[0, t] \,|\,B_t =0, L_t^0 = \ell)
	\stackrel{d}{=}
	(B[0, \tau_\ell] \,|\,\tau_\ell = t)
\end{equation}
for\/ $t = \lambda_j$ and\/
$\ell = \ell_j$.  \end{theorem}

\begin{proof}
Propositions~\ref{aldous:prpD},~\ref{aldous:prpTR}
and~\ref{aldous:pois} provide more explicit descriptions
of the law of $(\lambda_{I_j}, L^0_{I_j}, B^{\rm
br}_*[I_j])_{j =1,2, \ldots}$, for each of the four
interval partitions $(I_j)$. The above results for
occupation measures are deduced from these propositions
using Lemma~\ref{aldous:lmswap}. The fundamental {\em
switching identity}\/~\eqref{aldous:switch} is well known
\cite[\S 5]{aldous:py92}.  \qed\end{proof}

By general theory of local time processes
for diffusions or continuous semi-martingales
\cite{aldous:im65,aldous:rw87,aldous:ry99}, in the
Brownian and Bessel cases for each random subinterval $I$
of $[0,1]$ the random occupation measure $\mu_I$ derived
from $B^{\rm br}$ has an almost surely continuous density
$L_I^x$ relative to $m$ at $x$, where $m$ is a multiple
of the speed measure of the one-dimensional diffusion
$B$. To be precise about normalization of local times,
in the Brownian case with state space $\mathbb{R}$,
we take $m(\D x) = \D x$, so that~\eqref{aldous:taul}
holds with $\alpha = 1/2$ and $c = \sqrt{2}$.  In the
Bessel$(\delta)$ case with state space $\mathbb{R}_{\ge
0}$, we take $m(\D x) = 2 x^{\delta-1}\,\D x$, so
that~\eqref{aldous:taul} holds with $\alpha = 1 - \delta/2$
and $c = 2^{1 - \alpha} \Gamma(1-\alpha)/\Gamma(\alpha)$,
by \cite[(7.c)]{aldous:py92}.  In either case, $L^0_{I_j}$
in (iii) and (iv) is recovered like $\lambda_{I_j}$ as
a measurable function of the random occupation measure
$\mu_{I_j}$.  The distribution of the local time density of
the conditional occupation measure in (iv) is described
by a conditional form of the Ray--Knight theorem: see
\cite{aldous:leur96,aldous:jp97sde} for details in the
Brownian case.

\comment{This yields also the other
results stated in this section, in particular
Propositions~\ref{aldous:prpD},~\ref{aldous:prpTR}
and~\ref{aldous:pois} below which provide explicit
descriptions of the distributions of various functionals
of $B^{\rm br}$ derived from the $T$- and $D$-partitions.}

Our analysis of the $D$-partition is the following
expression of the decomposition of $B^{\rm br}$ at the
times $D_{V_j}$, implicit in \cite{aldous:ap92} in the
Brownian case:

\begin{lemma}[\cite{aldous:ap92}]\label{aldous:lmm-Drec}
For each\/ $j$, the pre-$D_{V_j}$ fragment of the bridge\/
$B^{\rm br}[0,D_{V_j}]$ is independent of the standardized
post-$D_{V_j}$ fragment\/ $B^{\rm br}_*[D_{V_j},1]$, which
has the same distribution  as\/ $B^{\rm br}$.  \end{lemma}

This is easily verified, because the $D_{V_j}$ are
stopping times relative to a filtration with respect to
which $B^{\rm br}$ has a strong Markov property.

\comment{\tt xxx well, maybe this would be tough to argue
in the general case. But one could do it then by excursion
theory anyway. leave as is?}

To describe various distributions, let $(\Gamma_s,
s \ge 0)$ denote a {\em standard gamma process},
that is the increasing L\'evy process with marginal
densities
\begin{equation}
\label{aldous:gams}
	P(\Gamma_s \in \D x)/\D x =  \frac{1}{\Gamma(s)}
	\,x^{s-1} \,\E^{-x} \qquad (x >0),
\end{equation}
so
$\Gamma_t - \Gamma_s \stackrel{d}{=} \Gamma_{t-s}$ for
$0 < s < t $.  Recall that for $a,b >0$ the beta$(a,b)$
distribution is that of
\begin{equation}
\label{aldous:bg}
	\mbox{$\beta_{a,b}:= \Gamma_a/\Gamma_{a+b}$, which is
	independent of $\Gamma_{a+b}$, with}
\end{equation}
\begin{equation}
\label{aldous:betab}
	P(\beta_{a,b} \in\D u)
	= \frac{\Gamma(a + b)}{\Gamma(a) \Gamma(b)}\, u^{a-1}
	(1-u)^{b-1}\,\D u  \qquad (0 < u < 1). 
\end{equation}
It is well known \cite[Lemma 3.7]{aldous:ppy92} that
for $G_t = G_t(B)$,
\begin{equation}
\label{aldous:dynk}
	\mbox{the standard $B$-bridge $B_*[0,G_t]$ is independent
	of $G_t$ with $G_t/t \stackrel{d}{=} \beta_{\alpha, 1
	- \alpha}$.}
\end{equation}
\comment{The next lemma is
\cite[Proposition 2]{aldous:ap92} in the Brownian case and
can be read from \cite[Proposition 15]{aldous:jp.bmpart}
in the Bessel case.}

\begin{lemma}[{\cite[Prop. 2]{aldous:ap92},
\cite[Prop. 15]{aldous:jp.bmpart}}]\label{aldous:lmmP2}
Let\/ $U$ with uniform$[0,1]$ distribution be independent
of\/ $B^{\rm br}$, and let\/ $G_U:= G_U(B^{\rm
br})$, $D_U:= D_U(B^{\rm br})$.  \comment{$$
I_1 :=
[0, G_U], \qquad I_2 := (G_U, D_U ], \qquad I_3 :=
(D_U,1]. 
$$} Then \comment{the random vector
of lengths\/ $(\lambda_{I_i}, i = 1,2,3)$ has the
Dirichlet\/ $(\alpha, 1-\alpha, \alpha)$ distribution,
that is}
$$
\comment{(\lambda_{I_1}, \lambda_{I_2},\lambda_{I_3})}
	(G_U, D_U - G_U, 1 - D_U)
	\stackrel{d}{=}
	(\Gamma_\alpha,\Gamma_{1} - \Gamma_{\alpha},
	\Gamma_{1 +\alpha} - \Gamma_1) /\Gamma_{1 + \alpha}. 
$$
Moreover,
the random vector\/ $(G_U, D_U - G_U, 1 - D_U)$ and the
three standardized processes\/ $B^{\rm br}_*[0,G_U]$,
$B^{\rm br}_*[G_U,D_U]$ and\/ $B^{\rm br}_*[D_U,1]$ are
independent, with
\begin{equation}
	B^{\rm br}_*[0,G_U]
	\stackrel{d}{=} B^{\rm br}_*[D_U,1] \stackrel{d}{=}
	B^{\rm br} \qquad\mbox{and}\qquad B^{\rm br}_*[G_U,D_U]
	\stackrel{d}{=} B^{\rm ex}. 
\end{equation}
\end{lemma}

\begin{proposition}\label{aldous:prpD} For the\/
$D$-partition

\noindent {\em (i)}\enspace the sequence of lengths is such
that
\begin{equation}
\label{aldous:llen}
	\lambda_{I^D_j} = W_j \prod_{i = 1}^{j-1} (1 - W_i)
\end{equation}
for a
sequence of independent random variables\/ $W_j$ with\/
$W_j \stackrel{d}{=} \beta_{1,\alpha}$.

\noindent {\em (ii)}\enspace The corresponding
sequence of local times at\/ $0$ can be expressed as
\begin{equation}
\label{aldous:llen1}
	L^0_{I^D_j} = \lambda_{I^D_j}^\alpha
	L_1^0\bigl(B^{\rm br}_*\bigl[ I^D_j\bigr]\bigr)
\end{equation}
where the\/ $L_1^0(B^{\rm br}_*[ I^D_j])$
are independent random variables, independent also of
the lengths\/ $\lambda_{I^D_j}$, with
\begin{equation}
\label{aldous:ltcon}
	L_1^0\bigl(B^{\rm br}_*\bigl[ I^D_j\bigr]\bigr)
\stackrel{d}{=} L_1^0(B) \stackrel{d}{=} \tau_1^{- \alpha}
\end{equation}
for\/ $\tau_1$ with the stable distribution
of index\/ $\alpha$ defined by \eqref{aldous:taul}.

\noindent {\em (iii)}\enspace The standardized path
fragments\/ $B^{\rm br}_*[ I^D_j]$ are independent
and identically distributed like\/ $B^{\rm
br}_*[0,D_U]$, and independent of the sequence of
lengths\/ $(\lambda_{I^D_j})$.  \comment{like\/ $B^{\rm
br}_*[0,D_U]$, where\/ $D_U:= D_U$ is the first zero of\/
$B^{\rm br}$ after an independent random time\/ $U$ with
uniform\/ $[0,1]$ distribution.}

\noindent {\em (iv)}\enspace For the length-ranked\/
$D$-intervals\/ $I^D_{(j)}$ instead of\/ $I^D_{j}$,
the lengths\/ $(\lambda_{I^D_{(j)}})$ have the
Poisson--Dirichlet$(\alpha)$ distribution defined by
ranking\/ $(\lambda_{I^D_j})$ as in\/ {\em (i)}, while
parts\/ {\em (ii)} and\/ {\em (iii)} hold without change.
\end{proposition}

\begin{proof} Parts (i)-(iii) are obtained by
repeated application of Lemmas~\ref{aldous:lmm-Drec}
and~\ref{aldous:lmmP2}, using the $\alpha$-scaling
rule~\eqref{aldous:locscale} for local times
and~\eqref{aldous:dynk}, as in \cite[Lemma
3.11]{aldous:ppy92}, for part (ii).  The second identity in
(ii) is a well-known consequence of the inverse relation
between $(L^0_t(B), t \ge 0)$ and $(\tau_\ell, \ell
\ge 0)$, as discussed in \cite{aldous:py92}.  Part (iv)
follows immediately from (i)-(iii).  \qed\end{proof}

See \cite{aldous:ki93,aldous:py95pd2} and
Lemma~\ref{aldous:poirep} for background on
the Poisson--Dirichlet distribution  appearing in
(iv).  L\'evy \cite{aldous:lev39} showed that in the
Brownian case the common distribution of $L_1^0(B)$ and
$\tau_1^{-1/2}$ appearing in~\eqref{aldous:ltcon} is simply
the distribution of $|B_1|$, with $B_1$ standard Gaussian.
But this does not generalize to the Bessel$(\delta)$
case for general $\delta = 2 - 2 \alpha$. Then $B_1
\stackrel{d}{=} \sqrt{2 \Gamma_{1 - \alpha}}$, which is a
simple transformation of the stable$(\alpha)$ distribution
of $\tau_1$ only for $\alpha = 1/2$.

The difficulty involved in Theorem~\ref{aldous:thmDT}
is that Definition~\ref{aldous:tdef} of the $T_j$
involves the local time $L^0_1:= L_1^0(B^{\rm
br})$, which depends on the path of $B^{\rm br}$
over the whole interval $[0,1]$. While we can
describe the finite-dimensional distributions
of the bivariate sequence $(\lambda_{I^T_{j}},
L^0_{I^T_{j}})_{j = 1,2, \ldots}$ by conditioning
on $L^0_1$ (see Proposition~\ref{aldous:prpT}),
this description is more complicated than our
description of $(\lambda_{I^D_{j}}, L^0_{I^D_{j}})_{j
= 1,2, \ldots}$ in Proposition~\ref{aldous:prpD}.
\comment{and it is not obvious how to pass from these
descriptions to the identity~\eqref{aldous:DTcore}.}
\comment{Lemma~\ref{aldous:lmm-Drec}, there is no simple
Markovian decomposition of the path of $B^{\rm br}$ at
time $T_j$, even for $j=1$.} In particular, \comment{Note
that~\eqref{aldous:DTcore} becomes false if $I^T_{(j)}$
replaced by $I^T_{j}$ and $I^D_{(j)}$ replaced by
$I^D_{j}$, because}
\begin{equation}
\label{aldous:neq}
	\lambda_{I^T_{1}} \stackrel{d}{\neq} \lambda_{I^D_{1}}.
\end{equation}
Indeed, by~\eqref{aldous:llen} we have
$$
	E\bigl(\lambda_{I^D_{1}}\bigr) = E(W_1) = 1/(1 + \alpha) > 1/2,
$$
whereas (by symmetry of $B^{\rm br}$ with respect to time
reversal in the Brownian or Bessel case) the distribution
of $T_1$ is symmetric about $1/2$, so whatever $\alpha
\in (0,1)$
$$
	E\bigl(\lambda_{I^T_{1}}\bigr) = E (T_1) = 1/2.
$$
Still, as explained combinatorially in the Brownian case
around~\eqref{aldous:locstick}, the two partitions give
rise to the same distribution for the sequence of local
times:
\begin{equation}
\label{aldous:Lsame}
	\Bigl(L^0_{I^D_j}\Bigr)
	\stackrel{d}{=}
	\Bigl(L^0_{I^T_j}\Bigr)
	:=
	\Biggl(L^0_{1}  \,U_j \prod_{i=1}^{j-1} (1 -U_i)\Biggr)
\end{equation}
where the second equality
by definition is read from~\eqref{aldous:tkdef}.
The first equality in distribution of sequences
follows from Lemma~\ref{aldous:lmm-Drec} and the
consequence of Lemma~\ref{aldous:lmmP2}, noted
in \cite[(3)-(4)]{aldous:ap92} in the Brownian
case, that
\begin{equation}
\label{aldous:Lsame1}
	\mbox{$L^0_{I^D_1}/L^0_1$ has uniform distribution on
	$(0,1)$, and is independent of $L^0_1$.}
\end{equation}
As indicated in Section~\ref{aldous:sec-IP}, this can
also be checked in general using the exchangeability
of the excursion interval partition.  \comment{On the
other hand this distribution can by represented (by
Theorem~\ref{aldous:thmDT}~(i)) as a $L^0$-biased sampling
from the sequence defined in~\eqref{aldous:llen},
and this representation is spelled out later as
Corollary~\ref{aldous:crlDT}.}

According to Proposition~\ref{aldous:prpD}, the
standardized bridge fragments over intervals of the
$D$-partition are i.i.d.\ copies of $B^{\rm br}_*[0,D_U]$,
both for the intervals in their original order and for
the intervals in length-ranked order.  A subtle feature
of the $T$-partition is that the standardized bridge
fragments over its intervals are neither independent nor
identically distributed in their original order, but these
fragments become i.i.d.\ when put into length-ranked
order.  This and other parallels between the $T$- and
$D$-partitions in length-ranked order are presented in
the following Proposition:

\begin{proposition}\label{aldous:prpTR} For the\/
$T$-partition in length-ranked order

\noindent {\em (i)}\enspace the sequence of
lengths\/ $(\lambda_{I^T_{(j)}})$ has the same
Poisson--Dirichlet$(\alpha)$ distribution as\/
$(\lambda_{I^D_{(j)}})$.

\noindent {\em (ii)}\enspace The corresponding sequence of
local times at\/ $0$ can be expressed as
\begin{equation}
\label{aldous:llen2}
	L^0_{I^T_{(j)}} =
	\lambda_{I^T_{(j)}}^\alpha \,
	L_1^0\bigl(B^{\rm br}_*\bigl[ I^T_{(j)}\bigr]\bigr)
\end{equation}
where the\/ $L_1^0(B^{\rm br}_*[ I^T_{(j)}])$ are
independent random variables, independent also of the
lengths\/ $L^0_{I^T_{(j)}}$, with
\begin{equation}
\label{aldous:ltcon1}
	L_1^0\bigl(B^{\rm br}_*\bigl[ I^T_{(j)}\bigr]\bigr)
	\stackrel{d}{=} L_1^0(B) \stackrel{d}{=} \tau_1^{- \alpha}
\end{equation}
just as in\/ \eqref{aldous:ltcon}.

\noindent {\em (iii)}\enspace The standardized
path fragments\/ $B^{\rm br}_*[ I^T_{(j)}]$ are
independent and identically distributed like\/ $B_* [0,
\tau_1]$, and independent of the sequence of lengths\/
$(\lambda_{I^T_{(j)}})$.  \end{proposition}

The only difference between this description of the
law of the sequence $(\lambda_{I_j}, L^0_{I_j}, B^{\rm
br}_*[I_j])_{j = 1,2, \ldots}$ for $I_j = I^T_{(j)}$, and
the previous description in Proposition~\ref{aldous:prpD}
for $I_j = I^D_{(j)}$, is that the common distribution
of the standardized $T$-fragments is that of $B^{\rm
br}_* [ 0, D_U]$, whereas the common distribution of
the standardized $D$-fragments is that of $B_* [0,
\tau_1 ]$.  The standardized process $B_* [0, \tau_1 ]$
is known as the {\em pseudo-bridge}\/ associated with the
self-similar Markov process $B$.  The following Lemma was
established by Biane, Le Gall and Yor \cite{aldous:bly87}
in the Brownian case, and extended to the Bessel case in
\cite[Theorem 5.3]{aldous:py92}.

\begin{lemma}[\cite{aldous:bly87,aldous:py92}]\label{aldous:lmmbly}
The law of the pseudo-bridge\/ $B_* [0, \tau_1 ]$ is
mutually absolutely continuous with respect to the law of
$B^{\rm br}$, with density proportional to\/ $1/L_1^{0}(B)$
relative to the law of\/ $B^{\rm br}$.  That is, for all
non-negative measurable path functionals\/ $F$
$$
	E \bigl[ F\bigl([B_* [0, \tau_1 ]\bigr) \bigr]
	= \frac{1}{c \alpha \Gamma(\alpha)}\,
	E \biggl[ \frac{F(B^{\rm br})}{L_1^{0} (B^{\rm br})}\biggr]. 
$$
where\/ $c$ is determined by the normalization
of local time via\/ \eqref{aldous:taul}.  \end{lemma}

While the laws of $B^{\rm br}_* [ 0, D_U]$ and
the pseudo-bridge $B_* [0, \tau_1 ]$ are mutually
singular,  their random occupation measures have the same
distributions.  In fact, the sample path of $B^{\rm br}_*
[ 0, D_U]$ is simply a random rearrangement of the sample
path of $B_* [0, \tau_1 ]$: \comment{Pairs of processes,
whose occupation density processes are identical
in law, are known as {\em local time equivalent}.
See \cite{aldous:xxx} for other interesting examples
of this phenomenon in the Brownian world, which have
been explained by similar path rearrangements, and
\cite{aldous:xxx} for closely related studies.}

\begin{lemma}\label{aldous:lmswap} Let\/ $U$ be a uniform\/
$(0,1)$ variable independent of\/ $B^{\rm br}$, and
independent of\/ $X$ distributed like\/ $B_* [0, \tau_1 ]$.
Then a process\/ $Y$ distributed like\/ $B^{\rm br}_* [
0, D_U]$ is created by the following rearrangement of the
path of\/ $X$, whereby the random occupation measures of\/
$X$ and\/ $Y$ are pathwise identical: let\/ $(G_U,D_U)$
be the excursion interval of\/ $X$ straddling time\/ $U$,
and let\/ $Y$ be derived from\/ $X$ by swapping the order
of the path fragments\/ $X[ G_U, D_U]$ and\/ $X[ D_U, 1]$,
say
\begin{equation}
\label{aldous:conc1}
	Y = X[0, G_U] : X[ D_U, 1] : X[G_U, D_U]
\end{equation}
with an obvious
notation for concatenation of path fragments.  \end{lemma}

\begin{proof} \comment{\paragraph{Proof of
Lemma~\ref{aldous:lmswap}}} By construction, the
path of $Y$ ends with a $B$-excursion of length $1 -
G_1(Y) = D_U - G_U$. The joint law of $Y[0,G_1(Y)]$
and $Y[G_1(Y),1]:= X[G_U,D_U]$ was described in
\cite[Theorem 1.3]{aldous:py92} and \cite[Theorem
3.1 and (3.d)]{aldous:ppy92}, and is identical to the
joint law of $Z[0,G_1(Z)]$ and $Z[G_1(Z),1]$ for $Z:=
B^{\rm br}_*[0,D_U(B^{\rm br})]$, which can be read from
Lemma~\ref{aldous:lmmP2}.  To be explicit, the common
distribution of $Y[0,G_1(Y)]$  and $Z[0, G_1(Z)]$ is
that of $B[0,G_1(B)]$ described by~\eqref{aldous:dynk},
\comment{, where $\smash{G_1(B) \stackrel{d}{=} \beta_{\alpha, 1 -
\alpha}}$,} while both $Y_*[G_1(Y),1]:= X_*[G_U,D_U]$ and
$Z_*[G_1(Y),1]:= B^{\rm br}_*[ G_U(B^{\rm br}),D_U(B^{\rm
br})]$ are standard $B$-excursions. Since the excursion
is in each case independent of the preceding fragment,
it follows that $\smash{Y \stackrel{d}{=} Z}$.  \qed\end{proof}

\comment{Part (i) of Theorem~\ref{aldous:thmDT}, and part
(iv) for the $D$- and $T$-partitions in their ranked
order, follow easily from Propositions~\ref{aldous:prpD}
and~\ref{aldous:prpTR}, using Lemma~\ref{aldous:lmswap}.
The remaining parts of Theorem~\ref{aldous:thmDT} are
then implied by the following alternate representation of
the law of the sequence of lengths and local times at $0$
for $D$- and $T$- partitions in their original order.}

\begin{proposition}\label{aldous:pois} Fix\/ $\xi >0$.
Let\/ $G$ be a random variable independent of\/ $B^{\rm
br}$, with\/ $G \stackrel{d}{=} \Gamma_\alpha/\xi$.
The distributions of the two bivariate sequences,
defined by the lengths and bridge local time measures of
intervals of the\/ $D$-partition and the\/ $T$-partition
respectively, are determined as follows:

\noindent {\em (i)}\enspace For\/ $I_j = I^D_j$ the
bivariate sequence
\begin{equation}
\label{aldous:poi1}
	\bigl(G  \lambda_{I_j},  G^\alpha  L^{\rm br}_{I_j}
	\bigr)_{j=1,2, \ldots}
\end{equation}
is the sequence of points\/
$(X_j,Y_j)$, in\/ $X$-biased random order, of a Poisson
process on\/ $\mathbb{R}_{>0}^2$ with intensity measure
\comment{at\/ $(x, y)$}
\begin{equation}
	\nu (\D t,\D\ell) :=
\comment{\frac12 x^{-1} \,\E^{- x/2}\,\D x\,P(\sqrt{x}\,|B_1| \in \D y)}
	\alpha t^{-1} \E^{- \xi t}\,\D t\,
	P (t^\alpha \tau_1^{- \alpha} \in \D\ell) =  \ell^{-1}
	P(\tau_\ell \in \D t) \E^{- \xi t}
\label{aldous:pp0}
\end{equation}
for\/ $\tau_1$ as in~\eqref{aldous:taul},
which makes
\begin{equation}
\label{aldous:gamdis}
	\Sigma_j X_j \stackrel{d}{=} \frac{\Gamma_\alpha}{\xi}
	\qquad\mbox{and}\qquad \Sigma_jY_j \stackrel{d}{=}
	\frac{\Gamma_1}{c  \xi^\alpha}. 
\end{equation}

\noindent {\em (ii)}\enspace If the points\/ $(X_j,Y_j)$
of a Poisson process with intensity\/ $\nu$ on\/
$\mathbb{R}_{>0}^2$ are listed in\/ $X$-biased order then
\begin{equation}
\label{aldous:normp}
	\bigl(\lambda_{I^D_j}, L^0_{I^D_j}\bigr)_{j = 1,2, \ldots}
	\stackrel{d}{=} \biggl(\frac{X_j}{\Sigma_X} ,
	\frac{Y_j}{\Sigma_X^\alpha} \biggr)_{j = 1,2, \ldots}
\end{equation}
for\/ $\Sigma_X:= \sum_j X_j$ as in\/ \eqref{aldous:gamdis}.

\noindent {\em (iii)}\enspace For\/ $I_j = I^T_j$
the bivariate sequence in\/ \eqref{aldous:poi1}
is the sequence of points, say\/ $(X_j',Y_j')$, in\/
$Y'$-biased random order, of another Poisson process
on\/ $\mathbb{R}_{>0}^2$ with the same intensity
measure\/ $\nu$.  So if in\/ {\em (ii)} the points\/
$(X_j,Y_j)$ are listed instead in\/ $Y$-biased order, then\/
\eqref{aldous:normp} holds with the sequence of\/
$T$-intervals instead of the sequence of\/ $D$-intervals.
\end{proposition}

\begin{proof} Part (i) is proved in
Section~\ref{aldous:sec-Dpart}.  Part (ii) is just
a restatement of part~(i).  Part (iii) is proved in
Section~\ref{aldous:sec-Tpart}.  \qed\end{proof}

Note that the normalization in~\eqref{aldous:normp}
involves $\Sigma_X$ and its $\alpha$th power, both for the
$D$-partition and for the $T$-partition. Obviously, this
is easier to handle if the sampling is $X$-biased rather
than $Y$-biased, which is one explanation of why various
distributions associated with $(I^D_j)$ are simpler than
their counterparts for $(I^T_j)$.

\section{Analysis of the $D$-partition}
\label{aldous:sec-Dpart}

As a preliminary for the proof of
Proposition~\ref{aldous:pois}~(i), we recall the following
well known lemma, which characterizes the distribution
of a sequence $(Q_j)$, known as the $GEM(\theta)$
distribution after Griffiths, Engen and McCloskey.
The distribution of $(Q_{(j)})$ obtained by ranking $(Q_j)$
is known as the {\em Poisson--Dirichlet distribution with
parameter\/ $\theta$}.  See \cite{aldous:ki75}, \cite[\S
9.6]{aldous:ki93}, \cite{aldous:py95pd2}.

\begin{lemma}[characterizations of\/ $GEM(\theta)$
\cite{aldous:mc65,aldous:ppy92,aldous:jp96bl}]\label{aldous:poirep}
Fix\/ $\theta >0$ and\/ $\xi > 0$.  Let\/ $G$ and\/
$Q_j, j = 1,2, \ldots$ be non-negative random variables.
Then the following are equivalent:

\noindent {\em (i)}\enspace the sequence\/ $(Q_j)$
admits the representation\/ $Q_{j} = W_{j} \prod_{i=
1}^{j-1} (1-W_{i})$ where the\/ $W_j$ are independent
beta$(1,\theta)$ variables, and\/ $G$ is independent of\/
$(Q_j)$ with\/ $\smash{G \stackrel{d}{=} \Gamma_\theta/ \xi}$;

\noindent {\em (ii)}\enspace $\sum_j Q_j = 1$ a.s. and\/
$(G Q_j)$ is the sequence of points of a Poisson point
process on\/ $\mathbb{R}_{>0}$ with intensity\/ $\theta
t^{-1} \E^{- \xi t}\,\D t$, listed in size-biased order.
\end{lemma}

The next well known result \cite[\S 5.2]{aldous:ki93},
\cite[Prop.  4.10.1]{aldous:resnickAdv}, combined with
the previous lemma, provides an efficient way to identify
various Poisson processes.

\begin{lemma}[Poisson marking]\label{aldous:poispray}
Let\/ $(S, \mathcal{S})$ and\/ $(T,\mathcal{T})$
be two measurable spaces.  Let\/ $(X_j)$ and\/ $(Y_j)$
be two sequences of random variables, with values in\/ $S$
and\/ $T$ respectively, such that the counting process\/
$\sum_j 1 (X_j \in \cdot)$ is Poisson with intensity
measure\/ $\mu$ on\/ $\mathcal{S}$, and the\/ $Y_j$ are
conditionally independent given\/ $(X_j)$, with
$$
	P(Y_j\in \cdot  \,|\,X_1, X_2, \ldots) =  P'(X_j, \cdot)
$$
for
some Markov kernel\/ $P'$ from\/ $(S, \mathcal{S})$ to\/
$(T, \mathcal{T})$.  Then the counting process\/ $\sum_j 1
((X_j,Y_j) \in \cdot)$ is a Poisson process on the product
space\/ $S\times T$ with intensity measure\/ $\mu(\D x)
P'(x,\D y)$ on the product\/ $\sigma$-field.  \end{lemma}

\begin{theopargself} \begin{proof}[\hskip-0.5em\
of Proposition~\ref{aldous:pois}
(i)] Proposition~\ref{aldous:prpD}~(i) and
Lemma~\ref{aldous:poirep}~(i) show that $(\lambda_{I^D_j},
j \geq 1)$ has $GEM(\alpha)$ distribution.  By assumption,
$G$ is independent of this sequence with $\smash{G \stackrel{d}{=}
\Gamma_\alpha/\xi}$.  Lemma~\ref{aldous:poirep} implies
that $(G \lambda_{I^D_j})$ is the size-biased ordering of
a Poisson point process of intensity $\alpha t^{-1}\E^{-
\xi t}\,\D t$.  Proposition~\ref{aldous:prpD}~(ii) and
Lemma~\ref{aldous:poispray} now identify the $(G
\lambda_{I^D_j},  G^\alpha  L^{\rm br}_{I^D_j})$
as the points of a Poisson process with intensity
measure $\nu$ defined by the first expression
in~\eqref{aldous:pp0}.  To check the equality of the two
expressions in~\eqref{aldous:pp0}, let
\begin{equation}
\label{aldous:fydef}
	f_\ell(t):= P(\tau_\ell \in \D t)/\D t.
\end{equation}
Since $\smash{\tau_\ell \stackrel{d}{=}
\ell^{1/\alpha} \tau_1}$ by~\eqref{aldous:taul},
\begin{equation}
\label{aldous:fy2}
	f_\ell(t) = \ell^{-1/\alpha} f_1(t/\ell^{1/\alpha})
\end{equation}
whereas by another change of variables
\begin{equation}
\label{aldous:chvar}
	P(t^\alpha \tau_1^{-\alpha}\in \D\ell)/\D\ell
	= \alpha^{-1} t \ell^{-1 - 1/\alpha}
	f_1(t/\ell^{1/\alpha}) = \alpha^{-1} t \ell^{-1} f_\ell(t)
\end{equation}
and the identity follows.  By application
of~\eqref{aldous:taul}, the $\ell$-marginal of $\nu$
is $\ell^{-1} \E^{- c \xi^\alpha \ell}\,\D\ell$.  The
distribution of $\sum_j Y_j$ is the infinitely divisible
law with this L\'evy measure, that is the exponential
distribution with rate $c \xi^\alpha$.  \qed\end{proof}
\end{theopargself}

Implicit in Lemma~\ref{aldous:lmmbly}
and~\eqref{aldous:chvar} is the following formula of
\cite[(3.u)]{aldous:ppy92} for the density of $L^0_1 :=
L^0_1 (B^{\rm br})$
\begin{equation}
\label{aldous:lbb}
	P(L^0_1 \in \D\ell) = c \alpha \Gamma(\alpha) \ell
	P(\tau_1^{-\alpha} \in \D\ell) = c \Gamma(\alpha)
	f_\ell(1) \,\D\ell
\end{equation}
for $f_\ell(x)$ as
in~\eqref{aldous:fydef} the stable$(\alpha)$ density of
$\tau_\ell$ determined by~\eqref{aldous:taul}.  That is
to say, the distribution of $L^0_1 (B^{\rm br})$ is
obtained by size-biasing the common distribution of
$\tau_1^{-\alpha}$ and $L^0_1(B)$.  In particular,
the general formula~\eqref{aldous:lbb} is consistent
with L\'evy's well known formulae in the Brownian
case \cite{aldous:lev39}, with $\alpha = 1/2,
c = \sqrt{2}$
\begin{equation}
\label{aldous:fellt}
	f_\ell(x) = \frac{\ell}{\sqrt{2 \pi}}\, x^{-3/2}
	\E^{-\frac12 \ell^2 / x}
\end{equation}
and
\begin{equation}
\label{aldous:fellt1}
	P(L^0_1 \in \D\ell)/\D\ell =
	\ell \E^{- \frac12 \ell^2}. 
\end{equation}
For general
$\alpha$, a series expression for $f_\ell(x)$ is known
\cite{aldous:po46,aldous:zo86,aldous:zolo94,aldous:uchzol99}.
If $\alpha = 1/n$ for some $n = 2,3, \ldots$,
integral expressions for $f_\ell(x)$ can be derived
from a representation of $1/\tau_\ell$ as a product
of $n-1$ independent gamma variables \cite[Theorem
3.4.3]{aldous:zo86}.  To conclude this section,
we record the following immediate consequence of
Proposition~\ref{aldous:pois}~(i):

\begin{corollary}\label{aldous:crlstab} Let\/
$(\sigma_j)$ be a sequence of i.i.d.\ copies of\/
$\tau_1$ with the stable\/ $(\alpha)$ distribution\/
\eqref{aldous:taul}, and let\/ $(Q_j)$ with\/
$GEM(\alpha)$ distribution of Lemma~\ref{aldous:poirep}
be independent of\/ $(\sigma_j)$.  Let\/ $L_j:=
(Q_j/\sigma_j)^\alpha$ and\/ $L:= \sum_j L_j$.  Then\/
$\smash{L \stackrel{d}{=} L^0_1}$ as in\/ \eqref{aldous:lbb},
and the sequence\/ $(L_j/L)$ has\/ $GEM(1)$ distribution,
independently of\/~$L$.  \end{corollary}

\section{Analysis of the $T$-partition}
\label{aldous:sec-Tpart}

We start by recalling the structure of a Markov process
up to the last time it visits its initial state before an
independent exponential time.  This does not involve the
self-similarity assumption.

\begin{lemma}[\cite{aldous:gp80}]\label{aldous:lmmgp}
Let\/ $(\tau_\ell, \ell \ge 0)$ be a drift-free
subordinator which is the inverse of the continuous local
time process\/ $(L^0_t(B), t \ge 0)$ of a regular recurrent
point\/ $0$, for a strong Markov process\/ $B$ started
at\/ $0$.  Let\/ $\varepsilon$ be an exponential variable
with rate\/ $\xi$, with\/ $\varepsilon$ independent of
$B$, and let
\begin{equation}
\label{aldous:GLdef}
	G :=G_\varepsilon (B) \qquad\mbox{and}\qquad L:= L_G^0(B)
	= L_\varepsilon^0(B). 
\end{equation}

\noindent {\rm (i)}\enspace The local time\/ $L$ has
exponential distribution with rate\/ $\psi(\xi)$, the
Laplace exponent of the subordinator defined by\/ $E(\E^{-
\xi \tau_\ell}) =  \E^{- \psi(\xi) \ell}$.

\noindent {\rm (ii)}\enspace For\/ $\ell > 0$, there
is the equality in distribution of path fragments
\begin{equation}
\label{aldous:chmeas}
	(B[0,G]  \,|\,L =\ell)
	\mathrel{\smash{\stackrel{d}{=}}}
	(B[0,\tau_\ell]  \,|\,\tau_\ell < \varepsilon). 
\end{equation}

\noindent {\rm (iii)}\enspace The joint distribution of\/
$(G,L)$ is
\begin{equation}
\label{aldous:taurho1}
	P(G\in \D t, L \in \D\ell)
	= \psi(\xi)\,\D\ell\, \E^{- \xi t}
	P( \tau_\ell \in \D t). 
\end{equation}
which is the
distribution of the value at time\/ $1$ of a drift free
bivariate subordinator with L\'evy measure
\begin{equation}
\label{aldous:taurho}
	\nu (\D t, \D\ell) = \ell^{-1}\,\D\ell\,
	\E^{- \xi t} P(\tau_\ell \in \D t)
\end{equation}
whose\/
$\ell$-marginal is the L\'evy measure\/ $\ell^{-1} \E^{-
\psi(\xi) \ell}\,\D\ell$ of the exponential distribution
of\/ $L$.  \end{lemma}

\begin{proof} These results are derived from It\^o's
theory of excursions of $B$, by letting $(N_t, t \ge 0)$
be a Poisson process with rate $\xi$, independent of $B$,
and taking $\varepsilon$ to be the time of the first point
of $N$. To briefly recall the argument, say that a jump
interval $(\tau_{y-}, \tau_y)$ of the inverse local time
process $\tau$ is {\em marked}\/ if $N(\tau_{y-}, \tau_y] >
0$ and {\em unmarked}\/ otherwise.  Then, by basic theory
of Poisson point processes, the sum of unmarked jumps
\begin{equation}
\label{aldous:subu}
	\tau^{\rm u}_\ell :=
	\sum_{0 < y < \ell} (\tau_{y} - \tau_{y-}) 1 (N(\tau_{y-},
	\tau_y] = 0)
\end{equation}
defines a subordinator
with distribution
\begin{equation}
\label{aldous:tauu}
	P(\tau^{\rm u}_\ell \in \D t) = \E^{\psi (\xi) \ell -
	\xi t} P(\tau_\ell \in \D t). 
\end{equation}
The left
end $G$ of the first marked interval is $G = \tau_{L -}
= \tau^{\rm u}_{L}$, and the subordinator $\tau^{\rm
u}$ summing unmarked jumps of $\tau$ is independent of
$L$, the local time of the first marked jump.  See also
\cite{aldous:gp80,aldous:py92,aldous:py97rh,aldous:rw87}.
\qed\end{proof}

To be more explicit, part (iii) of the Lemma states that
\begin{equation}
\label{aldous:gl}
	(G,L) \mathrel{\smash{\stackrel{d}{=}}}
	(\Sigma_j X_j, \Sigma_j Y_j)
\end{equation}
for $(X_j,Y_j)$ the points of a Poisson point process
on $\mathbb{R}_{>0}^2$ with intensity measure $\nu$
defined by~\eqref{aldous:taurho}.  In particular, for a
self-similar $B$ as in Section~\ref{aldous:sec-state},
this measure $\nu$ is identical to the measure $\nu$
featured in~\eqref{aldous:pp0}.

In the setting of Lemma~\ref{aldous:lmmgp}, even
with construction of the Poisson process of marks of
rate $\xi$ independent of $B$, more randomization
is required to construct points $(X_j,Y_j)$ such
that~\eqref{aldous:gl} holds with equality almost
surely rather than just in distribution.  But this can
be done by the following construction, which is the
basis of our proofs of Proposition~\ref{aldous:prpTR}
and Proposition~\ref{aldous:pois}~(iii).

\begin{lemma}\label{aldous:lmmgp2} In the setting of
Lemma~\ref{aldous:lmmgp}, let\/ $I:=  [G_I,D_I]$ be a
random subinterval of\/ $[0,1)$, where the endpoint\/ $1$
is deliberately excluded, to avoid the jump of the inverse
local time process\/ $(\tau_\ell)$ at time\/ $\ell = L$ in
the following construction.  Suppose\/ $I$ is independent
of\/ $B$ and\/ $\varepsilon$, and define further random
intervals
\begin{equation}
\label{aldous:indefs}
	I L :=
	[G_I L, D_I L ] \qquad\mbox{and}\qquad \tau(I L) :=
	[ \tau(G_I L), \tau(D_I L) ]
\end{equation}
where\/
$\tau(\ell):= \tau_\ell$ for\/ $\ell \ge 0$.

\noindent {\rm (i)}\enspace For\/ $y > 0$, there is the
equality in distribution of path fragments
\begin{equation}
\label{aldous:chmeas2}
	(B[ \tau(I L) ]  \,|\,\lambda_{IL}= y)
	\mathrel{\smash{\stackrel{d}{=}}}
	(B[0,\tau_y]  \,|\,\tau_y <
\varepsilon)
\end{equation}
where\/ $\lambda_{IL} :=
D_I L - G_I L$ is the increment of local time of\/ $B$
over the time interval\/ $\tau(I L)$.

\noindent {\rm (ii)}\enspace If\/ $(I_j)$ is an interval
partition of\/ $[0,1)$ which is independent of\/
$\varepsilon$ and\/ $B$, and\/ $(\tau_{I_j L})$  is the
corresponding interval partition of\/ $[0, G)$, then given
the sequence of local time increments\/ $(\lambda_{I_j
L})$ the path fragments\/ $B[ \tau_{I_j L} ]$ are
conditionally independent with distributions described
by\/ \eqref{aldous:chmeas2} for\/ $I_j$ instead
of\/ $I$.

\noindent {\rm (iii)}\enspace If\/ $I_j:=
[\widehat{V}_{j-1},\widehat{V}_j]$ with\/ $\widehat{V}_j:=
1 - \prod_{i = 1}^j(1-U_i)$ for independent uniform$(0,1)$
variables\/ $U_i$ independent of\/ $B$ and\/ $\varepsilon$,
then the bivariate sequence of local time increments
and path fragments
\begin{equation}
\label{aldous:ppp1}
	\bigl(\lambda_{I_j L}, B[ \tau(I_j, L) ]\bigr)_{j = 1,2, \ldots}
\end{equation}
is the sequence of points of a Poisson point
process on\/ $\mathbb{R}_{>0} \times \Omega$, in local-time
biased order, for a suitable space of path fragments\/
$\Omega$ of arbitrary finite length, with intensity
measure
\begin{equation}
\label{aldous:ppp2}
	\ell^{-1} P(\tau_\ell < \varepsilon, B[0, \tau_\ell ] \in \D\omega)
\end{equation}
whose\/ $\ell$-marginal is the L\'evy
measure\/ $\ell^{-1} \E^{- \psi(\xi) \ell}\,\D\ell$ of the
exponential distribution of\/ $L$ with rate\/ $\psi(\xi)$.

\noindent {\rm (iv)}\enspace Let\/ $X_j:= \tau(G_{I_j} L)
- \tau(G_{I_j} L)$ be the length and\/ $Y_j:= \lambda_{I_j
L}$ the local time increment associated with the random
subinterval\/ $\tau(I_j L)$ of\/ $[0,G)$. Then the\/
$(X_j,Y_j)$ are the points of a Poisson point process on\/
$\mathbb{R}_{>0}^2$ with intensity measure\/ $\nu$ defined
by\/ \eqref{aldous:taurho}, in\/ $Y$-biased random
order, and
\begin{equation}
\label{aldous:glas}
	(G,L) =(\Sigma_j X_j, \Sigma_j Y_j)
	\qquad\mbox{almost surely.}
\end{equation}
\end{lemma}

\begin{proof} The first two assertions are straightforward
consequences of the previous Lemma.  Part (iii) follows
from (ii), the Poisson representation of $GEM(1)$ in
Lemma~\ref{aldous:poirep}, and Poisson marking.  Part (iv)
follows from (iii) and Lemma~\ref{aldous:lmmgp}~(iii).
\qed\end{proof}

\begin{theopargself} \begin{proof}[\hskip-0.5em\ of
Proposition~\ref{aldous:pois}~(iii)] We will exploit the
following construction of the standard bridge $B^{\rm br}$
of the self-similar Markov process $B$ by random scaling,
as in \cite{aldous:py92} and \cite[Lemma 4]{aldous:py97rh}.
Let
\begin{equation}
\label{aldous:bcons}
	B^{\rm br}:= B_*[0, G_{\varepsilon}(B)]
	\qquad\mbox{for $\varepsilon$ independent of $B$ with}\qquad
	\varepsilon \mathrel{\smash{\stackrel{d}{=}}} {\Gamma_1 \xi},
\end{equation}
so $\varepsilon$ is exponential with rate
$\xi$. Then
\begin{equation}
\label{aldous:Gfacts}
	G :=
	G_\varepsilon (B) \stackrel{d}{=} \frac{\Gamma_\alpha}{\xi}
	\qquad\mbox{and}\qquad L:= L_G^0(B) = G^\alpha L^0_1(B^{\rm
	br}) \stackrel{d}{=} \frac{\Gamma_1}{c \xi^\alpha}
\end{equation}
by~\eqref{aldous:dynk},~\eqref{aldous:bg},
and $\alpha$-scaling of local times, where the exponential
distribution of $L$ is read from Lemma~\ref{aldous:lmmgp}.
Suppose now that $I_j:= [\widehat{V}_{j-1},\widehat{V}_j]$,
for $\widehat{V}_j$ as in Lemma~\ref{aldous:lmmgp2}.
The $T$-sequence is now constructed as a function
of these $\widehat{V}_j$ and $B^{\rm br}:=
B_*[0, G]$ as in~\eqref{aldous:bcons} according to
Definition~\ref{aldous:tdef}, that is
\begin{equation}
\label{aldous:trep}
	T_j := \inf \bigl\{u : L^0_u /LB_1 >\widehat{V}_j\bigr\}. 
\end{equation}
By~\ref{aldous:trep}
and~\eqref{aldous:Gfacts},
\begin{eqnarray*}
	\lambda_{I_j L}
	& = &
	\bigl(\widehat{V}_j - \widehat{V}_{j-1}\bigr) L =
	\bigl(\widehat{V}_j - \widehat{V}_{j-1}\bigr) \,G^\alpha L^0_1
	(B^{\rm br}) = G^\alpha L^0_{I^T_j}(B^{\rm br})
	\\
	\lambda_{\tau_{I_j L}}
	& = &
	\tau_{\widehat{V}_j L}
	- \tau_{\widehat{V}_{j-1} L} = G \,\lambda_{I^T_j}
	\\
	B_* [ \tau_{I_j L} ]
	& = &
	B^{\rm br}_*\bigl[ I^T_j \bigr].
\end{eqnarray*}
Part (iii) of Proposition~\ref{aldous:pois}
can now be read from Lemma~\ref{aldous:lmmgp2}~(iv).
\qed\end{proof} \end{theopargself}

\begin{theopargself} \begin{proof}[\hskip-0.5em\
of Proposition~\ref{aldous:prpTR}] Parts (i)
and (ii) follow immediately from the result of
Proposition~\ref{aldous:pois}~(iii) proved above.
Turning to consideration of the path fragments, we
observe by switching identity~\eqref{aldous:switch}
that
\begin{equation}
\label{aldous:switch1}
	(B_*[0,\tau_\ell ] \,|\,\tau_\ell = t)
	\mathrel{\smash{\stackrel{d}{=}}}
	(B^{\rm br} \,|\,L^0_1 = \ell t^{-\alpha})
\end{equation}
where
$L^0_1 := L^0_1(B^{\rm br})$ as usual, and a regular
conditional distribution for $B^{\rm br}$ given $L^0_1$
can be as constructed in \cite[Lemma 12]{aldous:jp97sde}.
Hence from~\eqref{aldous:ppp2}, if $\Omega_1$ denotes
a suitable space of paths of length 1, the trivariate
sequence of local time increments, lengths of path
fragments, and standardized path fragments
\begin{equation}
\label{aldous:seq1}
	\bigl(\lambda_{I_j L}, \lambda_{\tau_{I_j L}},
	B_* \bigl[ \tau_{I_j L} \bigr]\bigr)_{j = 1,2, \ldots}
	= \bigl(G^\alpha L^0_{I^T_j}(B^{\rm br}),
	G \lambda_{I^T_j}, B^{\rm br}_*\bigl[ I^T_j \bigr]
	\bigr)_{j = 1,2, \ldots}
\end{equation}
is a Poisson
process on $\mathbb{R}_{>0} \times \mathbb{R}_{>0} \times
\Omega_1$ whose intensity measure is
\begin{equation}
\label{aldous:seq2}
	\ell^{-1}\,\D\ell \, P(\tau_\ell \in \D t)
	\E^{- \xi  t} P(B^{\rm br} \in d \omega_1 \,|\,L^0_1 =
	\ell t^{-\alpha}). 
\end{equation}
Using the first form
of $\nu$ in~\eqref{aldous:pp0} to integrate out $\ell$
in~\eqref{aldous:seq2}, we see that the lengths and
standardized fragments
\begin{equation}
\label{aldous:seq3}
	\bigl(\lambda_{\tau_{I_j L}}, B_* [ \tau_{I_j L} ]
	\bigr)_{j = 1,2,\ldots}
	=\bigl(G \lambda_{I^T_j}, B^{\rm br}_*\bigl[ I^T_j \bigr]
	\bigr)_{j= 1,2, \ldots}
\end{equation}
form a Poisson process
on $\mathbb{R}_{>0} \times \Omega_1$ whose intensity
measure is
\begin{equation}
\label{aldous:pp3}
	\alpha t^{-1} \E^{- \xi t}\,\D t\,Q(\D\omega_1)
\end{equation}
where
\begin{equation}
\label{aldous:pp4}
	Q(\D\omega_1)
	= \int_0^\infty P(B^{\rm br} \in \D\omega_1 \,|\,L^0_1
	= y) P ( \tau_1^{-\alpha}  \in \D y)
	= P(B_*[0,\tau_1] \in \D\omega)
\end{equation}
by the switching
identity~\eqref{aldous:switch1}.  The factorization
in~\eqref{aldous:pp3} shows that the $B^{\rm br}_*[ I^T_j
]$ are i.i.d.\ copies of $B_*[0, \tau_1]$ when listed in
length-ranked order. That is part~(iii) of Proposition
\ref{aldous:prpTR}.  \qed\end{proof} \end{theopargself}

\subsection{Further distributional results}

We  record in this section a number of further formulae
related to the distribution of the lengths and local times
defined by the $T$-partition.

\begin{proposition}\label{aldous:prpT} For the\/
$T$-partition the\/ $(2n+1)$-variate joint density of
the total bridge local time\/ $L^0_1$, the lengths of
the first\/ $n$ intervals, and the local times at\/ $0$
on these intervals, is given by the formula
$$
\displaylines{\indent
	P\bigl(L^0_1\in \D\ell, \lambda_{I^T_j} \in \D x_j,L^{\rm br} (I^T_j)
	\in d y_j, 1 \le j \le n\bigr)
	\hfill\cr\hfill{}
	=
	c \Gamma(\alpha)\,\D\ell
	\,f_{\ell - y_1 - \cdots - y_n}(1-x_1- \cdots
	-x_n) \prod_{j = 1}^n \frac{\D x_i \,\D y_j\, f_{y_j}(x_j)}
	{\ell - y_1 - \cdots - y_{j-1}}
	\indent}
$$
for\/ $f_y(x):=
P(\tau_y \in \D x)/\D x$ the stable$(\alpha)$ density as
in\/ \eqref{aldous:fydef}.  \end{proposition}

\begin{proof} This follows from the switching
identity~\eqref{aldous:switch} and the definition of the
$T$-partition.  \qed\end{proof}

While the distributions of the cut times $T_k$
and interval lengths $(T_k - T_{k-1})$ in principle
determined Proposition~\ref{aldous:prpT}, formulae for
these distributions are more easily obtained as follows.
For $0 < u < 1$, let
$$
	\tau^{\rm br}_u  := \inf \{t : L^0_t /L^0_1 = u \}
$$
where $(L^0_t, 0 \le t \le 1)$ is
the local time process at $0$ of $B^{\rm br}$.  Then by
use of the switching identity~\eqref{aldous:switch}
we can write down for $0 < x <1$, $0 < \ell < \infty$,
\begin{equation}
\label{aldous:taubu}
	P(\tau^{\rm br}_u \in \D x \,|\,L^0_1 = \ell) /\D x
	=
\comment{f(x \,|\,u,\ell) :=}
	\frac{f_{u \ell}(x) f_{\bar{u} \ell} (\bar{x})}{f_{\ell}(1)}
\end{equation}
where $\bar{x} :=
1 - x$.  \comment{For each fixed $u$, this function
of $x$ integrates to $1$ over $0 \le x\le 1$ by the
convolution rule $f_v * f_w = f_{v+w}$.} Integrating out
with respect to the distribution~\eqref{aldous:lbb}
of $L^0_1$ gives the density
\begin{equation}
\label{aldous:taubdens}
	P(\tau^{\rm br}_u \in \D x)/\D x
\comment{= f(x \,|\,u) :=}
	= c \Gamma(\alpha)\int_0^\infty  f_{u \ell}(x)
	f_{\bar{u} \ell} (\bar{x})\,\D\ell. 
\end{equation}
which can be simplified using
L\'evy's formula~\eqref{aldous:fellt} in the Brownian
case to give for $\alpha = 1/2$
\begin{equation}
\label{aldous:fxu}
\comment{f(x \,|\,u)}
	P(\tau^{\rm br}_u \in \D x) /\D x
	= \frac{u \,\bar{u}}{2 (\bar{x} u^2 +x \bar{u}^2)^{3/2}}
	\qquad (0 < u, x < 1). 
\end{equation}
\comment{As a check, this integrates to $1$ over $0 <
x <1$ by elementary calculus.} In particular, for $u
= 1/2$ we recover the the result of \cite[Theorem
3.2]{aldous:bp92} that $\tau^{\rm br}_{1/2}$ has uniform
distribution on $[0,1]$ in the Brownian case.

According to Definition~\ref{aldous:tdef}, $T_k :=
\tau^{\rm br}_{\widehat{V}_k}$, for $\widehat{V}_k$
independent of $B^{\rm br}$ with
$$
	1 - \widehat{V}_k
	\mathrel{\smash{\stackrel{d}{=}}}
	\widehat{V}_k - \widehat{V}_{k-1}
	\mathrel{\smash{\stackrel{d}{=}}} \Pi_k
$$
for $\Pi_k$ is a product of $k$
independent uniform$(0,1)$ variables, with
\begin{equation}
\label{aldous:ppdens}
	\frac{P(\Pi_k  \in \D u)}{\D u} =
	\frac{(- \log u)^{k-1}}{(k - 1)!}
	\qquad\mbox{and}\qquad
	\sum_{k=1}^\infty \frac{P(\Pi_k \in \D u)}{\D u}
	= \frac{1}{u}
\end{equation}
because $\log \Pi_k $
is the $k$'th point of a rate $1$ Poisson process on
$[0,\infty)$.  Since the process $(\tau^{\rm br}_u, 0 \le
u \le 1)$ has exchangeable increments, we find that $1 -
T_k$ and the length of the $k$th $T$-interval have the
common distribution
\begin{equation}
\label{aldous:tk2}
	P(1 - T_k  \in \D x)
	= P\bigl(\lambda_{I^T_k} \in \D x\bigr)
	=\int_0^1 P(\tau^{\rm br}_u \in \D x) P(\Pi_k \in \D u).
\end{equation}
In particular,  in the Brownian case $\alpha
= 1/2$,~\eqref{aldous:tk2} and~\eqref{aldous:fxu} yield
the curious formula
\begin{equation}
\label{aldous:taud}
	\frac{P(T_1 \in \D x)}{\D x}
	=\frac{h(x) + h (\bar{x})}{2}
	\qquad\mbox{with}\qquad
	h(x):= \frac{1}{\sqrt{x}} + \log
	\biggl(\frac{1}{\sqrt{x}} - 1  \biggr).
\comment{\frac{1}{2}
\left[ \frac{1}{\sqrt{x}} + \frac{1}{\sqrt{\bar{x}}} +
\log \left(\frac{\sqrt{x} - x}{\sqrt{\bar{x}} + \bar{x}}
\right) \right]}
\end{equation}

\begin{corollary}\label{aldous:crlspec} The point
process of lengths of\/ $T$-intervals has mean density
\begin{equation}
\label{aldous:specdens}
	\sum_{k = 1}^\infty P\bigl(\lambda_{I^T_k} \in \D x\bigr)
	=  \alpha  x^{-1}(1-x)^{\alpha - 1}\,\D x
	= \int_0^1 P(\tau^{\rm br}_u\in \D x) u^{-1}\,\D u
\end{equation}
for\/ $x \in (0,1)$.
\end{corollary}

\begin{proof} The first equality is read from part
(i) of Proposition~\ref{aldous:prpTR} and the well
known formula for the mean density of points of a
Poisson--Dirichlet$(\alpha)$ distributed sequence
\cite[(6)]{aldous:py95pd2}, which can be read from
Lemma~\ref{aldous:poirep}.  The second equality is then
read from~\eqref{aldous:tk2} and~\eqref{aldous:ppdens}.
\qed\end{proof}

For general $\alpha$,
the second equality in~\eqref{aldous:specdens}
does not seem very obvious from~\eqref{aldous:taubu}
and~\eqref{aldous:taubdens}.  However, it can be checked
for $\alpha = 1/2$ using \eqref{aldous:fxu}, and it
can also be verified by a very general argument, which we
indicate in Section~\ref{aldous:sec.intensity}.

Path decompositions of $B^{\rm br}$ at the times
$T_k$ are more complicated than the corresponding
decompositions for the times $D_{V_j}$ expressed by
Lemma~\ref{aldous:lmm-Drec}.  For the $T$-partition,
the pieces are not pure $B$-bridges. Rather, when
normalized they have density factors involving their
local times at $0$.  Compare with similar constructions in
\cite{aldous:bly87,aldous:fpy92,aldous:ppy92,aldous:py97kt}.

By the Poisson analysis of the previous section,
conditionally given $(T_1, \allowbreak L^0_{T_1}, L^0_1)$
the pieces of $B^{\rm br}$ before and after time $T_1$
are independent $B$-bridges with prescribed lengths and
local times at $0$.  The appearance of $h+k$ in formula
(a) below shows that the right side does not factor into
a function of $(x,h)$ and a function of $(x,k)$.  So even
in the Brownian case, $L^0_{T_1}$ and $L^0_1 - L^0_{T_1}$
are not conditionally independent given $T_1$, and hence
the same can be said of the fragments of $B^{\rm br}$
before and after time $T_1$.

\begin{proposition}\label{aldous:PdrT2} In the Brownian
case with\/ $\alpha = 1/2$, $c = \sqrt{2}$,
$$
	\frac{P(T_1 \in \D x, L^0_{T_1} \in \D h, L^0_1 - L^0_{T_1}
	\in \D k)}{\D x \,\D h \,\D k} = \frac{h \,k}{\sqrt{2
	\pi}}\,\frac{(x \bar x)^{-3/2}}{(h + k)} \exp \biggl(-
	\frac{h^2}{2x} - \frac{k^2}{\bar{x}} \biggr)
$$
while
for
$$
	X:= \frac{L^0_{T_1}}{\sqrt{T_1}} = L_1^0 (B^{\rm br}_* [ 0,T_1])
	\qquad\mbox{and}\qquad
	Y:= \frac{L^0_{1} - L^0_{T_1}}{\sqrt{1 - T_1}}
	= L_1^0 (B^{\rm br}_* [ T_1,1]). 
$$
there is the joint density
\begin{equation}
\label{aldous:fab}
	\frac{P(X \in \D a, Y \in \D b)}{\D a \,\D b}
	= \frac{a \,b}{\sqrt{2 \pi}}\, I(a,b)
	\exp \biggl(- \frac{a^2}{2} - \frac{b^2}{2} \biggr)
\end{equation}
where
$$
	I(a,b):= \int_0^1
	\frac{(x \bar x)^{-1/2}}{a \sqrt{x}  + b \sqrt{\bar{x}}}\,\D x
	= \frac{1}{r} \log \biggl(  \frac{(r + a) (r + b)}{(r -
	a) (r - b)} \biggr)
$$
for\/ $ r:= \sqrt{a^2 + b^2} $.
\end{proposition}

\begin{proof} The first formula is an instance of
Proposition~\ref{aldous:prpT} which we now check.
With notation as in~\eqref{aldous:taubu},
$$
	P(T_1 \in\D x \,|\,L^0_{T_1} = h, L^0_1  - L^0_{T_1} = k)/\D x
	= f(x \,|\,u, \ell)
$$
for $h = u \ell$ and $ k = \bar{u}
\ell$.  We also know, by definition of $T_1$, that
$$
	P (L^0_{T_1} \in \D h, L^0_1 - L^0_{T_1} \in \D k) =
	\D h \,\D k \, \E^{-\frac12 \ell^2}
$$
where $\ell =
h + k$ and an $\ell^{-1}$ has canceled the factor of
$\ell$ in the density~\eqref{aldous:fellt1} of $L^0_1$.
Combining these formulae gives the trivariate density of
$(T_1, L^0_{T_1}, L^0_1 - L^0_{T_1})$, which rescales to
give
$$
	\frac{P(T_1 \in \D x, X \in \D a, Y  \in \D b)}
	{\D x \,\D a \,\D b}
	= \frac{a \,b}{\sqrt{2 \pi}}\,
	\frac{(x \bar x)^{-1/2}}{(a \sqrt{x}  + b \sqrt{\bar{x}})}
	\exp \biggl(- \frac{a^2}{2} - \frac{b^2}{2} \biggr).
$$
and~\eqref{aldous:fab} follows by integrating out $x$.
\qed\end{proof}

\section{Complements} \label{aldous:sec-COMP}

\subsection{Mappings conditioned to have a single basin}

In the Brownian case, a variation of the transformation
from $X$ to $Y$ in Lemma~\ref{aldous:lmswap},
which further swaps the exchangeable pair of
fragments $X[0, G_U]$ and $X[ D_U, 1]$, is the
continuous analog of the transformation, mentioned in
Fact~\eqref{aldous:comfacts}(e) from the stretch of
the cycles-first mapping walk for a given basin to the
stretch of the basins-first walk for the same basin.
As pointed out in the last section of \cite{aldous:ap92},
if the uniform mapping of $[n]$ is conditioned to have
only one cycle, the scaled basins-first walk converges in
distribution to the process $2 |B^{\rm br}_*[0,D_U(B^{\rm
br})]|$.  The above argument yields:

\begin{corollary}\label{aldous:crlmax} For a uniform
mapping of\/ $[n]$ conditioned to have only one cycle,
the scaled cycles-first walk converges in distribution
to\/ $ 2 |B_* [ 0, \tau_1]|$ where\/ $B_* [ 0, \tau_1]$
is the Brownian pseudo-bridge.  \end{corollary}

The distributions of several basic functionals of
pseudo-bridge $B_* [ 0, \tau_1]$ are known.  In particular,
the occupation density of the reflected process is governed
by the same stochastic differential equation governing the
occupation density  process of a reflecting Brownian bridge
or Brownian excursion \cite{aldous:jp97sde}.  According to
Knight \cite{aldous:kni88} (see also \cite{aldous:py97kt}
and papers cited there), the law of the maximum of the
reflected pseudo-bridge is identical to that of $1/(2
\sqrt{H_1(R_3)})$ where $H_1(R_3)$ is the hitting time
of $1$ by the three-dimensional Bessel process, with
transform $E(\exp(- \frac12 \theta^2 H_1(R_3))) = \theta/
\sinh \theta$ for real $\theta$.  Thus we deduce:

\begin{corollary}\label{aldous:crlmax1} For a uniform
mapping of\/ $[n]$ conditioned to have only one cycle,
the asymptotic distribution of the maximum height of any
tree above the cycle, normalized by\/ $\sqrt{n}$, is the
distribution of\/ $1/\sqrt{H_1(R_3)}$.  \end{corollary}

See also \cite{aldous:bpy99z} for a survey of closely
related distributions and their applications.

\subsection{Exchangeable interval partitions}
\label{aldous:sec-IP}

Suppose that $(I^{\rm ex}_j)$ is an exchangeable interval
partition of $[0,1]$.  That is (assuming for simplicity
that the lengths $\lambda_{I^{\rm ex}_j}$ are almost
surely all distinct), for each $n = 2,3, \ldots$ such
that $\lambda_{I^{\rm ex}_{(n)}} > 0$, where $(I^{\rm
ex}_{(j)})$ is the associated length-ranked interval
partition, conditionally given $\lambda_{I^{\rm ex}_{(n)}}
> 0$ the ordering of the longest $n$ sub-intervals $I^{\rm
ex}_{(j)}, 1 \le j \le n$ is equally likely to be any
one of the $n!$ possible orders, independently of the
lengths of these $n$ intervals.  Call $(I^{\rm ex}_j)$
\comment{{\em finite}\/ if $P(\lambda_{I^{\rm ex}_{(n)}}
> 0) \to 0$ as $n \to \infty$, in which case we denote by
$J$ the almost surely finite number of intervals in the
partition, and call $(I^{\rm ex}_j)$} {\em infinite}\/ if
$P(\lambda_{I^{\rm ex}_{(n)}} > 0) = 1$ for all $n$.  As
shown by Kallenberg \cite{aldous:kal83l}, for an infinite
exchangeable interval partition $(I^{\rm ex}_j)$, for each
$u \in [0,1]$ the fraction of the longest $n$ intervals
that lie to the left of $u$ has an almost sure limit
$\bar{L}^0_u$ as $n \to \infty$. The process $(\bar{L}^0_u,
0 \le u \le 1)$ is a continuous increasing process,
the {\em normalized local time process of\/ $(I^{\rm
ex}_j)$}. It is easily shown that for $B^{\rm br}$ as in
previous sections, and more generally for $B^{\rm br}$
the standard bridge of any nice recurrent Markov process,
constructed as in \cite{aldous:fpy92}, the interval
partition $(I^{\rm ex}_j)$ defined by the excursions of
$B^{\rm br}$ away from $0$ is an infinite exchangeable
interval partition of $[0,1]$, whose normalized local time
process is $\bar{L}^0_u = L^0_u/L^0_1, 0 \le u \le 1$ for
any of the usual Markovian definitions of a bridge local
time process $L^0_u := L^0_u(B^{\rm br})$. In particular,
this remark applies to a self-similar recurrent process $B$
as considered in previous sections.

\begin{theorem}\label{aldous:thmex} The assertions of
Theorem~\ref{aldous:crl1} remain valid for the\/ $D$- and\/
$T$-partitions defined by Definitions~\ref{aldous:ddef}
and~\ref{aldous:tdef} for any infinite exchangeable
interval partition\/ $(I^{\rm ex}_j)$ instead of the
excursion intervals of a standard Brownian bridge\/ $B^{\rm
br}$, with the complement of\/ $\bigcup_j I^{\rm ex}_j$ in\/
$[0,1]$ instead of the zero set of\/ $B^{\rm br}$, and
the normalized local time process\/ $(\bar{L}^0_u, 0 \le u
\le 1)$ of\/ $(I^{\rm ex}_j)$ instead of\/ $(L^0_u/L^0_1,
0 \le u \le 1)$.  Moreover, the sequence of normalized
local times\/ $(\bar{L}^0_{I_j})$ has the same\/ $GEM(1)$
distribution for\/ $I_j = I^D_j$ as for\/ $I_j = I^T_j$.
\end{theorem}

Theorem~\ref{aldous:thmex} can be derived from a
certain combinatorial analog, stated and proved as
Lemma~\ref{aldous:lmmex} below.  Let us briefly outline
the method of derivation, without details.  Consider an
infinite exchangeable interval partition $(I_j)$.
Take $k$ independent uniform $(0,1)$ sample points,
assign ``weight" $1/k$ to each, and let $(I^{(k)}_j)$
be the intervals containing at least one sample point.
Each interval $I^{(k)}_j$ is thereby assigned weight
$1/k \times $ (number of sample points in interval).
For fixed $k$ we can apply Lemma~\ref{aldous:lmmex},
interpreting ``length" as ``weight", and conditionally on
the number of intervals in the partition.  The conclusion
of Lemma~\ref{aldous:lmmex} is a variant of the desired
Theorem~\ref{aldous:crl1} for $(I_j)$, in which ``position
$x \in (0,1)$" of interval endpoint is replaced by
``$1/k \times $ (number of sample points in $(0,x)$)",
and in which ``normalized local time at $u \in (0,1)$"
is replaced by ``relative number of sampled intervals
in $(0,u)$".  One can now argue that as $k \to \infty$
we have a.s. convergence of these variant quantities to
the original quantities in Theorem~\ref{aldous:thmex}.

\begin{lemma}\label{aldous:lmmex} Let\/ $(I^{\rm
ex}_i)_{1 \le i \le n}$ be an exchangeable interval
partition of\/ $[0,1]$ into\/ $n$ subintervals of
strictly positive length. Define\/ $D_{V_j}$ as in
Definition~\ref{aldous:ddef} for\/ $1 \le j \le J^D_n$,
where\/ $J_n$ is the first\/ $j$ such that\/ $D_{V_j} = 1$,
to create a\/ $D$-partition\/ $(I^D_j)_{1 \le j \le J^D_n}$
of\/ $[0,1]$, and define a\/ $T$-partition\/ $(I^T_j)_{1
\le j \le J^T_n}$ of\/ $[0,1]$ similarly using cut points\/
$T_j, 1 \le j \le J^T_n$ determined as follows: given that
the random set of endpoints of\/ $(I^{\rm ex}_j)_{1 \le
j \le n}$ is\/ $\mathcal{U}:= \{u_j\}_{0 \le j \le n}$
with\/ $0 = u_0 < u_1 < \cdots < u_n = 1$, let\/ $T_1$
have uniform distribution on\/ $\mathcal{U} \cap (0,1]$,
and given also\/ $T_1 = t_1 < 1$ let\/ $T_2$ have uniform
distribution on\/ $\mathcal{U} \cap (t_1, 1]$, and so
on, until\/ $T_{J^T_n} = 1$. For\/ $I_j$ an interval of
either of the\/ $D$- or\/ $T$-partitions so defined, let\/
$N_{I_j}$ denote the number of intervals of\/ $(I^{\rm
ex}_i)_{1 \le i \le n}$ which are contained in\/ $I_j$,
so\/ $1 \le N_{I_j} \le n$.  Then the assertions of
Theorem~\ref{aldous:crl1} remain valid provided that\/
$N_{I_j}$ is substituted everywhere for $L^0_{I_j}$.
\end{lemma}

\begin{proof} We will check that part (i) of
Theorem~\ref{aldous:crl1} holds in this setup, along
with~\eqref{aldous:perms}. The remaining claims are
straightforward and left to the reader.  By conditioning
on the ranked lengths $\lambda_{I^{\rm ex}_{(j)}}$ of the
intervals $(I^{\rm ex}_i)_{1 \le i \le n}$, it suffices to
consider the case when these ranked lengths are distinct
constants.  Let $\Pi_n^D$ denote the random partition of
$[n]$ defined by the random equivalence relation $i \sim
j$ iff $I^{\rm ex}_{(i)}$ and $I^{\rm ex}_{(j)}$ are part
of the same component interval of the $D$-partition, and
define $\Pi_n^T$ similarly in terms of the $T$-partition.
Since each unordered collection of lengths and sub-interval
counts is a function of the corresponding partition, it
suffices to show that $\Pi_n^D \stackrel{d}{=} \Pi_n^T$.
Due to the well known connection between the discrete
stick-breaking scheme used to define the $T$-partition
and the cycle structure of random permutations,
which was recalled in Section \ref{aldous:comfacts}
(d), we can write down the distribution of $\Pi_n^T$
without calculation: for each unordered partition of
$[n]$ into $k$ non-empty subsets $\{A_1, \ldots, A_k
\}$,
\begin{equation}
\label{aldous:PTform}
	P\bigl(\Pi_n^T =\{A_1, \ldots, A_k \}\bigr)
	= \frac{1}{n!} \prod_{j = 1}^k\bigl(|A_j| - 1\bigr) ! 
\end{equation}
where $|A_i|$ is the
number of elements of $A_i$.  On the other hand, for the
$D$-partition, for each ordered partition $(A_1, \ldots,
A_k)$ and each choice of $a_j \in A_j$,$ 1 \le j \le
k$, with $\lambda(a)$ the length of $I^{\rm ex}_{(a)}$
and $\lambda(A) := \sum_{a \in A} \lambda(a)$, we can
write down the probability
$$
	P\biggl(I^D_j = \!\!\bigcup_{a \in A_j}\!\!
	I^{\rm ex}_{(a)}  \ \mbox{and}\ I^{\rm ex}_{(a_j)} \
	\mbox{has right end}\ D_{V_j} \biggr)
	= \frac{1}{n!} \prod_{j= 1}^k \bigl(|A_j| - 1\bigr)!\,
	\frac{\lambda(a_j)}{\sum_{i = j}^k \lambda(A_i)}
$$
where the factors of $(|A_j| - 1) ! $
come from the different possible orderings of all but
the last $I^{\rm ex}_{(i)}$ to form $I^D_j$.  If we now
sum over all choices of $a_j \in A_j$, for each $1 \le j
\le k$, we find that $\lambda(a_j)$ is simply replaced by
$\lambda(A_j)$. If we then replace $(A_1, \ldots, A_k)$ by
$(A_{\sigma(1)}, \ldots, A_{(\sigma(k)})$ and sum over all
permutations $\sigma$ of $[k]$, to consider all sequences
of sets consistent with a given unordered partition $\{A_1,
\ldots, A_k \}$, we get precisely~\eqref{aldous:PTform}
for $\Pi_n^D$ instead of $\Pi_n^T$, due to the identity
$$
	\sum_{\sigma} \prod_{j = 1}^k
	\frac{\lambda(A_{\sigma(j)})}{\sum_{i = j}^k \lambda (A_{\sigma(i)})}
	= 1. 
$$
This is
obvious, because the product is the probability of picking
the sequence of sets $(A_{\sigma(j)}, 1 \le j \le k)$ in
a process of $\lambda(A_i)$-biased sampling of blocks of
the partition $\{A_1, \ldots, A_k \}$.  \qed\end{proof}

We note the consequence of the previous proof that the number
of components $J^D_n$ of the $D$-partition and the number
of components $J^T_n$ of the $T$-partition have the same
distribution, which is the same for every exchangeable
interval partition $(I^{\rm ex}_i)_{1 \le i \le n}$ of
$[0,1]$ into $n$ subintervals of strictly positive length:
\begin{equation}
\label{aldous:perms}
	J^D_n \stackrel{d}{=}
	J^T_n \stackrel{d}{=} K_n \stackrel{d}{=}
	\sum_{i = 1}^n 1_{C_i}
\end{equation}
where $K_n$ is the number of
cycles of a uniformly distributed random permutation of
$[n]$, and the $C_i$ are independent events with $P(C_i)
= 1/i$.  The second two of these equalities in distribution
are well known and easily explained without calculation
\cite{aldous:jp96bl}.  But the first is quite surprising,
and we do not see how to explain it any more simply than
by the previous proof.

\subsection{Intensity measures}
\label{aldous:sec.intensity}

In this section we check Corollary~\ref{aldous:crlspec}
by showing it can be generalized and proved as follows:

\begin{corollary}\label{aldous:crlspec1} In the setting of
Theorem~\ref{aldous:thmex}, the common intensity measure
of the the point process of lengths of\/ $T$-intervals
and the the point process of lengths of\/ $D$-intervals
is
\begin{equation}
\label{aldous:specdens1}
	\sum_{k= 1}^\infty P\bigl(\lambda_{I^T_k} \in \D x\bigr)
	= \sum_{k =1}^\infty P\bigl(\lambda_{I^D_k} \in \D x\bigr)
	= \frac{P (D_{V_1}\in \D x)}{x}
	= \int_0^1 \frac{P(\bar{\tau}_u \in \D x)}{u}\,\D u
\end{equation}
where\/ $\bar{\tau}_u:= \inf
\{t : \bar{L}^0_t > u \}$ is the inverse of the normalized
local time process of the exchangeable interval partition.
\end{corollary}

\begin{proof} The equality of the first three measures
displayed in~\eqref{aldous:specdens1} is read from the
conclusion of Theorem~\ref{aldous:thmex}, using the
fact that the $D$-partition is in length-biased order.
The equality of the first and fourth measures follows from
the definition of the $T_k$, the exchangeable increments of
$(\bar{\tau}_u, 0 \le u \le 1)$, and~\eqref{aldous:ppdens},
just as in the proof of~\eqref{aldous:specdens}.
\qed\end{proof}

As a check on Theorem~\ref{aldous:thmex},
let us verify the equality of the second and
fourth measures in~\eqref{aldous:specdens1} in the
following special case, which includes the setting of
Corollary~\ref{aldous:crlspec}.

Let $(\tau_\ell, \ell \ge 0)$ be the inverse
local time process of $B$ at $0$, for $B$ as in
Lemma~\ref{aldous:lmmgp} not necessarily self-similar.
Note that we must explicitly assume $(\tau_\ell, \ell \ge
0)$ is drift free for the conclusion of part (iii) of that
Lemma to be true.  We assume that now.  Assume that the
L\'evy measure of $(\tau_\ell, \ell \ge 0)$ has density
$\rho(x)$.  Let $(I^{\rm ex}_j)$ be the exchangeable
partition of $[0,1]$ generated by the excursion intervals
of $B$ conditional on $B_1 = 0$ and $L_1(B) = \ell$ for
some fixed $\ell > 0$, or equivalently by the jumps of
$(\tau_s, 0 \le s \le \ell)$ given $\tau_\ell = 1$.  Then,
formula~\eqref{aldous:taubu} generalizes easily to show
that the fourth measure in \eqref{aldous:specdens1} has
density at $x$
\begin{equation}
\label{aldous:thirdmeas}
	\int_0^1 u^{-1}\,\D u\,
	\frac{f_{u \ell}(x) f_{\bar{u}
	\ell} (\bar{x})}{f_{\ell}(1)}
\end{equation}
for
$f_\ell(x)$ as in~\eqref{aldous:fydef}.  On the other hand,
abbreviating $D:= D_{V_1}$ and $G:= G_{V_1}$ so $[G,D]$ is
the interval $I^{\rm ex}_j$ which covers the independent
uniform time $V_1$, we know from~\eqref{aldous:two} that
for $0 < w < 1$
$$
	P\bigl(1 - (D-G) \in \D w\bigr)
	= \frac{\ell \rho(1 - w) (1-w) f_\ell(w)\,\D w}{f_\ell(1)}\,. 
$$
Also, it is
easily seen that conditionally given $ 1 - (D-G) = w$, the
normalized local time $\bar{L}^0_G$ is uniform on $(0,1)$
and independent of the pair $(G,1-D)$, which is distributed
like $(\tau_{u\ell}, \tau_{\bar{u} \ell})$ conditioned
on $\tau_\ell = w$.  Together with the previous formula
for $w = y + 1-x$, this gives the trivariate density
$$
	\frac{P\bigl(\bar{L}^0_G \in \D u , G \in \D y, D \in \D x\bigr)}
	{\D u \,\D y \,\D x}
	= \frac{\ell \rho(x - y) (x-y) f_{u\ell}(y)
	f_{\bar{u} \ell}(\bar{x})}{f_\ell(1)}
$$
($0 < y < x <1$).
Now~\eqref{aldous:two} implies that
$$
	\int_0^y
	f_{u\ell}(y) \rho(x - y) (x-y)\,\D y =
	\frac{x f_{u\ell}(x)}{u \ell}
$$
so we deduce that
$$
	\frac{P\bigl(\bar{L}^0_G \in
	\D u , D \in \D x\bigr)}{\D u  \,\D x}
	= \frac{x f_{u\ell}(y)
	f_{\bar{u} \ell}(\bar{x})}{u f_\ell(1)}
$$
and hence
that the density displayed in~\eqref{aldous:thirdmeas}
is indeed $x^{-1}P(D \in \D x)/\D x$.

\subsection{Two orderings of a bivariate Poisson process}

According to Proposition~\ref{aldous:pois}, for each
$\alpha \in (0,1)$ the Poisson point process with
intensity measure $\nu(\D x,\D y) = \rho(x,y)\,\D x\,\D
y$ displayed in~\eqref{aldous:pp0} has the following
paradoxical property:	\begin{description}[(a)]

\item[(a)] If the
points $(X_j,Y_j)$ are put in $X$-biased order, then the
$Y_j$ are in $Y$-biased order, whereas

\item[(b)] if the points $(X_j,Y_j)$ are put in $Y$-biased
order, then the $X_j$ are not in $X$-biased random order;
even the distribution of $X_1$ is wrong.

\end{description} We first see this for $\alpha = 1/2$
by passage to the limit of elementary combinatorial
properties of uniform random mappings.  We then see it
for general $\alpha$ from the bridge representations of
Proposition~\ref{aldous:pois}.  Other point processes
of lengths and local times with these properties can
be constructed from an exchangeable interval partition,
as shown by Theorem~\ref{aldous:thmex} in the previous
section and Lemma~\ref{aldous:poirep}.  This argument,
shows that (a) holds for the bivariate Poisson process
with intensity \eqref{aldous:taurho} featured in
Lemma~\ref{aldous:lmmgp}, for any drift free subordinator
$(\tau_y,  y \ge 0)$ with $E(\E^{- \xi \tau_y}) =  \E^{-
\psi(\xi) y}$.  Then the $Y_j$ normalized by their sum
have $GEM(1)$ distribution, both for an $X$-biased and for
a $Y$-biased ordering.  We offer here a slightly different
explanation of (a) in this case. That is, given some joint
density $\rho(x,y)$, we indicate conditions on $\rho$
which are necessary and sufficient for (a) to hold for
the bivariate Poisson process with intensity $\rho$, and
then check that these conditions are in fact satisfied in
the case~\eqref{aldous:taurho}.

Let $(X_j,Y_j)$ be the points of a Poisson process on
$\mathbb{R}_{>0}^2$ with intensity $\rho(x,y)\,\D x\,\D y$,
in $X$-biased order.  Let $\Sigma_X : = \sum_j X_j$ and
$\Sigma_Y : = \sum_j Y_j$.  Let
$$
	f_X (x) := P(\Sigma_X\in \D x) /\D x ;
	\qquad f_Y (y) := P(\Sigma_Y \in \D y)/\D y
$$
$$
	\rho_X(x):= \int_{0}^\infty \rho(x,y)\,\D y;
	\qquad \rho_Y(y):= \int_{0}^\infty \rho(x,y)\,\D x. 
$$
By a basic Palm calculation, as in \cite{aldous:ppy92}
\begin{equation}
\label{aldous:two}
	P(X_1 \in \D x,\Sigma_X - X_1 \in \D w)
	= \rho_X(x)\,\D x \,f_X(w)\,\D w\,\frac{x}{x + w}
\end{equation}
and similarly, with
$$
	f_{X,Y} (x,y) := P(\Sigma_X \in \D x, \Sigma_Y \in \D y)
	/ (\D x\,\D y)
$$
\begin{multline}
\label{aldous:three}
	P(X_1 \in \D x, Y_1 \in \D y, \Sigma_X - X_1 \in \D w,
	\Sigma_Y - Y_1 \in \D v)
	\\
	= \rho(x,y) \,\D
	x \,\D y \,f_{X,Y}(w,v) \,\D w \,\D v \,\frac{x}{x + w}\,.
\end{multline}
Now, a {\em necessary condition}\/ for the $Y_j$ to be in
$Y$-biased order is that $Y_1$ should have the same joint
distribution with $\Sigma_Y$ as if $Y_1$ were a size-biased
pick from the $Y_i$, that is like~\eqref{aldous:two}
\begin{equation}
\label{aldous:four}
	P(Y_1 \in \D y,\Sigma_Y - Y_1 \in \D v)
	= \rho_Y(y)\,\D y \,f_Y(v)\,\D v \,\frac{y}{y + v}\,. 
\end{equation}
Thus a necessary
condition on $\rho(x,y)$ for (a) to hold is that for
all $y,v \ge 0$
\begin{equation}
\label{aldous:five}
	\int_0^\infty \D x \int_0^\infty \D w\, \rho(x,y)
	f_{X,Y}(w,v) \,\frac{x}{x + w} = \rho_Y(y) f_Y(v)
	\,\frac{y}{y + v}. 
\end{equation}
Moreover, by keeping track of
the first $k$ of the $(X_j,Y_j)$ jointly with $\Sigma_X$
and $\Sigma_Y$ it is clear that we can write down a
multivariate version of \eqref{aldous:five} whose truth
for all $k$ would be necessary and sufficient for (a).

In the special case~\eqref{aldous:taurho}, with $f_y(x):=
P(\tau_y \in \D x)/\D x$, the subordination argument gives
$$
	\rho(x,y) = y^{-1} f_y(x) \E^{- \xi x}. 
$$
Since the
$Y$-marginal is exponential with rate $\psi(\xi)$,
$$
	f_Y(y) = \psi(\xi)  y \rho_Y(y) = \psi(\xi)
	\E^{- \psi (\xi)  y}
$$
and hence by generalization
of~\eqref{aldous:chmeas}, using $\smash{(\Sigma_X, \Sigma_Y)
\stackrel{d}{=} (G,L)}$,
$$
	f_{X,Y}(x,y) = \psi(\xi)
	f_y(x) \E^{- \xi x}  = \psi(\xi) y \rho(x,y). 
$$
If these
expressions are substituted in~\eqref{aldous:five}, and
we use the definition of $\psi(\xi)$ on the right side,
we find that~\eqref{aldous:five} reduces to the identity
$$
	E \biggl[\frac{\tau_y}{{\tau_y + \tau_v} \E^{- \xi
	(\tau_y + \tau_v)}} \biggr]
	= \frac{y}{y + v}\, E \bigl[
	\E^{- \xi (\tau_y + \tau_v)} \bigr]. 
$$
But this is true
by virtue of
$$
	E \biggl[ \frac{\tau_y}{{\tau_y +
	\tau_v}} \biggm| \tau_y + \tau_v \biggr]
	= \frac{y}{y + v}
$$
which holds by exchangeability of increments of
$(\tau_\ell, \ell \ge 0)$.  Moreover, the multivariate
form of~\eqref{aldous:five} mentioned above is easily
checked the same way.

\begin{acknowledgement} We thank Gregory Miermont and
an anonymous referee for careful reading and helpful
comments.\end{acknowledgement}

\end{document}